\documentclass[12pt,reqno]{amsart}
\usepackage{latexsym}
\usepackage{amsmath,amsthm,amsfonts,amscd,eucal,amssymb}
\usepackage{epsfig}
\usepackage{enumerate}
\usepackage{xcolor}
\usepackage{xcolor}
\usepackage[all]{xy}
\usepackage{pdfsync}

\usepackage{tikz}

\newtheorem{thm}{Theorem}[section]
\newtheorem{cor}[thm]{Corollary}

\newtheorem{lem}[thm]{Lemma}

\newtheorem{prop}[thm]{Proposition}

\theoremstyle{definition}
\newtheorem{dfn}[thm]{Definition}

\theoremstyle{remark}
\newtheorem{rem}[thm]{Remark}
\newtheorem{ex}[thm]{Example}
\numberwithin{equation}{section}
\newtheorem*{ack}{Acknowledgments}

\def\ca{{\mathcal A}}

\def\cc{{\mathcal C}}

\def\cp{{\mathcal P}}

\def\car{{\mathcal R}}

\def\ct{{\mathcal T}}

\def\bc{{\mathbb C}}

\newcommand{\bn}{\mathbb N}
\def\N{{\mathbb N}}
\def\br{{\mathbb R}}

\renewcommand{\a}{\alpha}

        \def\G{\Gamma}

\def\s{\sigma}

\def\f{\varphi}

        \def\O{\Omega}

\def\ov{\overline}

\newcommand{\hooklongrightarrow}{\lhook\joinrel\longrightarrow}

\begin{document}

\title{Self-similar topological fractals}
\author{Fabio E.G. Cipriani}
\address{(F.E.G.C.) Politecnico di Milano, Dipartimento di Matematica,
piazza Leonardo da Vinci 32, 20133 Milano, Italy.} \email{fabio.cipriani@polimi.it}
\author{Daniele Guido}
\address{(D.G.) Dipartimento di Matematica, Universit\`a di Roma ``Tor
Vergata'', I--00133 Roma, Italy.} \email{guido@mat.uniroma2.it}
\author{Tommaso Isola}
\address{(T.I.) Dipartimento di Matematica, Universit\`a di Roma ``Tor
Vergata'', I--00133 Roma, Italy} \email{isola@mat.uniroma2.it}
\author{Jean-Luc sauvageot}
\address{(J.-L.S.) Institut de Math\'ematiques de Jussieu - Paris rive gauche, CNRS -  Universit\'e Paris Cit\'e - Universit\'e Paris Sorbonne, F 75205 PARIS Cedex 13, France.}
\email{jean-luc.sauvageot@imj-prg.fr}


\date{}

\begin{abstract}
We introduce the notion of (abelian) similarity scheme, as a constructive model for topological self-similar fractals, in the same way in which the notion of iterated function system furnishes a constructive notion of self-similar fractals in a metric environment. At the same time, our notion gives a constructive approach to the Kigami-Kameyama notion of topological fractals, since a similarity scheme produces a topological fractal a la Kigami-Kameyama, and many Kigami-Kameyama topological fractals may be constructed via  similarity schemes. Our scheme consists of objects $X_0\stackrel{\f}{\hooklongrightarrow}X_1\stackrel{\pi}{\leftarrow} Y\times X_0$, where $X_0,X_1$ and $Y$ are compact Hausdorff spaces, the map $\f$ is continuous injective and the map $\pi$ is continuous surjective. This scheme produces a sequence $X_n$, $n\in\bn$, of compact Hausdorff spaces, $X_n$ embedded in $X_{n+1}$, and a compact Hausdorff space $X_\infty$ giving a sort of injective limit space, which turns out to be self-similar. We observe that the space $Y$ parametrizes the generalized similarity maps, and  finiteness of $Y$ is not required.
\end{abstract}
\maketitle

\section{Introduction}
While the notion of fractal has been long popularized, so that many laypeople may associate to such word either broccoli or the coast of Norway, a generally accepted precise notion of fractal has never been achieved. Often, fractals are characterised in the negative: their dimension is not integer, they are not manifolds {\it a la } Withney, since, due to ramification points, they are not locally homeomorphic to an Euclidean ball, they cannot be studied via the notion of $G$-space {\it a la } Klein, since their simmetry group, either topologically or metrically, is usually finite.

However, things change dramatically when self-similar fractals are concerned: their description as a space which is locally homeomorphic to itself brought to the precise notion of IFS (iterated function system) due to Hutchinson, where the fractal is given by means of a complete metric space (the ambient space) and a finite number of contraction endomorphisms which give rise to a compact subspace $K$ (the fractal) of the ambient space. In this setting, the key role is played by the ambient space and the metric. Such notion has many advantages: it is constructive, namely the few data mentioned above produce a metric fractal space, whose properties are not immediately evident from the outset, and allows various finitary approximations, with which an analytic study of the fractal can be performed.

The aim of this paper is to devise an analogous notion for topological self-similar fractals, where neither a metric nor an ambient space is required, but still few data produce a fractal together with a family of approximating spaces. 

Let us remark here that the notion of topological self-similar fractal (or set) has been considered by Kigami and Kameyama \cite{Ka,Ki}, where such fractal is described as a quotient of the infinite address space $Y^\bn$, $Y=\{1,2,\dots, n\}$, with the natural shift maps on the address space descending to continuous (injective) maps on the quotient. This notion is very general, indeed contains IFS fractals, but is not constructive, the fractal set being part of the definition.

\subsubsection*{The classical notion of self-similar fractal {\it a la} Hutchinson}
The standard notion of self-similar fractal consists of a complete metric space $(\Omega,d)$ and of a finite set $(f_1,\dots f_n)$ of contraction maps. Then one obtains a map $\Phi$ on the subsets of $\Omega$ given by $\Phi(C)=\cup_i f_i(C)$. Clearly $\Phi$ on compact subsets is a contraction w.r.t. the Hausdorff metric, and the fractal $K$ associated with $\Phi$ is simply its unique fixed point:
\begin{equation}\label{fixed}
\Phi(K)=K.
\end{equation}
We remark that, by the Banach fixed point theorem, given any compact subset $X_0\subset \Omega$,
\begin{equation}\label{limit}
\lim_n \Phi^n(X_0)=K,
\end{equation}
where the limit is taken w.r.t. the Hausdorff metric. For some pcf (post-critically finite) fractals, $X_0$ may be taken as the set of essential fixed points, thus obtaining an increasing sequence $X_n=\Phi^n(X_0)$, which approximates the fractal in the Hausdorff metric.

We note that the metric on $\Omega$ plays here a crucial role.

\subsubsection*{What if we remove the metric?}

We may then consider a locally compact space   $\Omega$  together with continuous maps $(f_1,\dots f_n)$,  $f_j:\Omega\to\Omega$. Assuming the existence of a compact subspace $X_0\subset \Omega$  such that $X_1:=\Phi(X_0)\supset X_0$, we obtain an increasing family $\{X_n=\Phi(X_{n-1})\}$ of compact sets. Then the fractal $K$ could be defined as the closure of $\cup_n X_n$  in  $\Omega$.
Of course, $K$ is not necessarily compact in general, however, if we assume so, $K$ is self-similar, i.e. $\Phi(K)=K$. 
Indeed, setting $A=\cup_n X_n$, we get $\Phi(A)=A$. Therefore $\Phi(K)$ is a compact, hence closed, set containing $\Phi(A)=A$, therefore it contains $K$. Conversely, by the continuity of the maps $\{ f_i,i=1,\dots,n\}$, $\Phi(\overline A)\subset \overline{\Phi( A)}$, and $K$ is a fixed point.
In particular,  when there exists a metric for which the maps $f_i$ are contractions, the two definitions give the same object.
We note that in both definitions, the existence of an ambient space appears crucial.

\subsubsection*{Do we really need an ambient space?}

The first place where $\Omega$ plays a role is in the definition of $\Phi$, namely in the way in which the various copies $\{f_i(X_0),i=1,\dots n\}$ of $X_0$ embed in $\Omega$. The lack of an ambient space can be resolved by specifying a gluing procedure, which means a projection map from the disjoint union $\sqcup_i f_i(X_0)$ to $X_1$. This amounts to consider the diagram
\begin{equation}
\begin{matrix}
& &Y\times X_0
\\
& &\downarrow\pi\\
X_0&\stackrel{\f}{\hooklongrightarrow}&X_1,
\end{matrix}
\end{equation}
where $X_0,X_1,Y$ are compact Hausdorff spaces (indeed $Y$ was the set $\{1,\dots,n\}$ up to now), $\f$ is continuous injective, and $\pi$ is continuous surjective. Such a diagram will be called an {\it (abelian) similarity scheme}.

Such a scheme corresponds to the map $\Phi$ of the previous descriptions, namely it allows to associate to a pair $(Z,\f_Z)$, where $Z$ is compact Hausdorff and $\f_Z:X_0\to Z$ is continuous injective, a new pair 
$(\widehat{Z},\f_{\widehat{Z}})=\Phi(Z,\f_{Z})$ according to the commutative diagram
\begin{equation}\label{FixedPoint}
\begin{matrix}
& &Y\times X_0&\stackrel{ id_Y \times \f_Z }{\hooklongrightarrow}&Y\times Z
\\
& &\downarrow\pi& &\downarrow\pi_Z\\
X_0&\stackrel{\f}{\hooklongrightarrow}&X_1&\stackrel{\f_{\widehat{Z},X_1}}{\hooklongrightarrow}&\widehat{Z}
\end{matrix}
\end{equation}
where $\pi_Z$ is the quotient map w.r.t. the minimal equivalence relation $\sim$ on $Y\times Z$ such that $\pi(y, x_0)=\pi(y', x'_0)\Rightarrow(y, \f_Z(x_0))\sim (y', \f_Z(x'_0))$. In this way, the map $\f_{\widehat{Z},X_1}$ is well defined, and we set $\f_{\widehat{Z}}= \f_{\widehat{Z},X_1} \circ\f$. 

\begin{rem}
Let us observe that we no longer require $Y$ to be a finite set, our analysis indeed works for any compact Hausdorff space $Y$.
\end{rem}

\subsubsection*{Self-similar fractals as fixed points of an endo-functor}

It is then natural to define a self-similar set associated with the similarity scheme $X_0\stackrel{\f}{\hookrightarrow} X_1\stackrel{\pi}{\leftarrow}Y\times X_0$ as a fixed point for the map $\Phi$, namely as a pair $(Z,\f_Z)$ which is isomorphic to 
$\Phi(Z,\f_{Z})$.
From the technical point of view, the map from $(Z,\f_Z)$ to $(\widehat{Z},\f_{\widehat{Z}})$ is an endofunctor of the category $\ct_{X_0}$ whose objects are pairs $(Z,\f_Z)$, with $Z$ compact Hausdoff and $\f_Z:X_0\to Z$ continuous injective, and a morphism $f_\vartheta:(Z,\f_Z)\to(W,\f_W)$ is given by a continuous injective map $\vartheta:Z\to W$ such that $\f_W=\vartheta\circ\f_Z$. It is easy to see that for pcf fractals with $X_0,X_1,Y$ as above, the self-similar fractal is indeed a fixed point according to the definition above.

However, up to now, we haven't given a procedure to obtain a fixed point associated with a similarity scheme, and it is time to face this problem.

\subsubsection*{An approximation procedure}
Mimicking the limit construction, cf. \eqref{limit}, we may apply the scheme to the pair $(X_1,\f)$, obtaining the pair $(X_2,\f_{2,0})=\Phi(X_1,\f)$. A
repeated application of this procedure gives rise to an infinite diagram
\begin{equation}
\hskip-0.8cm \begin{matrix}
&                                                                 
&Y\times X_0&\stackrel{id_Y \times \f_{1,0}}{\hooklongrightarrow}
&Y\times X_1&\stackrel{id_Y \times \f_{2,1}}{\hooklongrightarrow}
&Y\times X_2&\stackrel{id_Y \times \f_{3,2}}{\hooklongrightarrow}
&Y\times X_3&\stackrel{id_Y \times \f_{4,3}}{\hooklongrightarrow}
\\
& 
&\downarrow\pi_0& 
&\downarrow\pi_1& 
&\downarrow\pi_2& 
&\downarrow\pi_3& 
\\
   X_0&\stackrel{\f_{1,0}}{\hooklongrightarrow}
&X_1&\stackrel{\f_{2,1}}{\hooklongrightarrow}
&X_2&\stackrel{\f_{3,2}}{\hooklongrightarrow}
&X_3&\stackrel{\f_{4,3}}{\hooklongrightarrow}
&X_4&\stackrel{\f_{5,4}}{\hooklongrightarrow}
\end{matrix}
\end{equation}
where all the squares are commutative, all vertical maps are surjective, and all horizontal maps  are injective.

We gained what looks like an approximation procedure, therefore we may expect to get a candidate of a fixed point as a kind of limit of the $X_n$ above. However, this is the second point where the lack of an ambient space is a problem: we cannot set $X_\infty$ as the closure of the union of the $X_n$'s. Instead, we would like to define $X_\infty$ as a compactification of the union $\cup X_n$, but it is not clear which compactification to choose. In other terms, the challenge now is to define a suitable compact space $X_\infty$ in which all the $X_n$ consistently embed via a family of coherent injective maps, namely $X_\infty$ should satisfy: 

\begin{equation}\label{embeddings}
\begin{matrix}
\exists\f_{\infty,n}:X_n\to X_\infty \textrm{ such that}\\
\f_{\infty,n}=\f_{\infty,m}\circ\f_{m,n}\\
X_\infty=\overline{\cup_n\f_{\infty,n}(X_n)}.
\end{matrix}
\end{equation}
\subsubsection*{Ideas from Kigami-Kameyama notion of topological self-similar fractal}
According to Kameyama, \cite{Ka}, see also \cite{Ki}, the fractal $K$ is  a quotient of the space $Y^\bn$ of infinite addresses endowed with the product topology, with $Y=\{1,\dots,n\}$, endowed with shift maps $f_y:K\to K$, $y\in Y$, such that the following diagram commutes:
\begin{equation}\label{shifts}
\begin{matrix}
Y^\bn&\stackrel{W_y}{\longrightarrow}&Y^\bn\\
\downarrow&&\downarrow\\
K&\stackrel{f_y}{\longrightarrow}& K,
\end{matrix}
\end{equation}
where $W_y(y_1,\dots ,y_n,\dots)=(y,y_1,\dots, y_n,\dots)$.

Following this suggestion, we look for $X_\infty$ as a quotient of $Y^\bn$ w.r.t. a suitable (closed) equivalence relation, endowed with the quotient topology.

Under the mild hypotesis of full injectivity (see below), we define such an equivalence relation, prove that the $X_n$'s and $X_\infty$ satisfy Property \ref{embeddings} and  prove that such limit point is indeed a fixed point of the endo-functor. In this way, even though uniqueness for fixed points does not hold in general (cf. Example \ref{nonunique}), we may associate with any similarity scheme a canonically defined limit fixed point $X_\infty=\pi^\bn(Y^\bn)$, therefore a similarity scheme provides a constructive approach to a topological fractal (the limit fixed point) together with a family of approximating spaces $X_n$. These results constitute the main part of this research, cf. Theorem \ref{LimFixPt}, and are contained in Sections 3, 4, 5. 

In particular, Section 3 contains a definition of a combinatorial equivalence relation ${\mathcal R}_\infty$ such that $Y^\bn/{\mathcal R}_\infty$ is a limit space of the sequence $X_n$, namely satisfies Property \ref{embeddings}, cf. Theorem \ref{embeddings2}.
\\
Section 4 shows that, assuming full injectivity, the relation ${\mathcal R}_\infty$ is closed, so that $X_\infty$ is compact Hausdorff, cf. Theorem \ref{closedRel}.
\\
Section 5 shows that the limit space $X_\infty$ is indeed a fixed point of the endo-functor, namely it is self-similar, cf. Theorem \ref{LimFixPt}.

In Section 6 we  compare the topological limit fractals produced by similarity schemes with the topological fractals obeying the Kigami-Kameyama requests. 
On the one hand, for any such limit fractal we get a commutative diagram
\begin{equation}
\begin{matrix}
Y^\bn&\stackrel{W_y}{\longrightarrow}&Y^\bn\\
\downarrow\pi^\bn&&\downarrow\pi^\bn\\
X_\infty&\stackrel{f_y}{\longrightarrow}& X_\infty,
\end{matrix}
\end{equation}
so that, when $Y$ is finite, we obtain a Kigami-Kameyama topological fractal, cf. Corollary \ref{LimFixPtIsKKfractal}.
On the other hand, given a Kigami-Kameyama topological fractal $K$, it is possible to  choose $X_0\subset K$ such that $(X_0, X_1=\cup_{j=1}^n f_j(X_0),\{1,\dots,n\})$ is a similarity scheme, and its limit fixed point, together with its similarity maps,  coincides with the given topological fractal, cf. Theorem \ref{thm:KK}.

Finally, various examples of similarity schemes together with their associated limit fixed points are described in Section 7.

We conclude this introduction with various observations.

We first remark that our interest for a constructive notion of a topological fractal independent of a metric and of an ambient space is also due to our search of a notion of noncommutative self-similar fractal (or self-similar $C^*$-algebra). Indeed the notion of similarity scheme can easily be dualized and transported to the noncommutative: a noncommutative (or $C^*$) similarity scheme is given by a diagram 
\begin{equation}
\begin{matrix}
& &B\otimes A_0
\\
& &\uparrow
\\
A_0&\stackrel{\f}{\longleftarrow}&A_1
\end{matrix}
\end{equation}
where $A_0,A_1,B$ are unital $C^*$-algebras, $\f$ is a surjective morphism and the vertical arrow is an inclusion. All the constructions that we give below for an abelian similarity scheme may be performed in the non-abelian setting, giving rise to noncommutative fractals as suitable limit fixed points. An example of such construction can be recognized in the paper \cite{CGIS3}, where the noncommutative gasket is constructed from the scheme
\begin{equation}
\begin{matrix}
& &M_3(\bc)\otimes \bc^3
\\
& &\uparrow
\\
\bc^3&\stackrel{\f}{\longleftarrow}&\bc^6.
\end{matrix}
\end{equation}
The description of noncommutative self-similar fractals together with their properties is contained in the work in progress \cite{CGIS5}.

We then note that the idea of constructing self-similar fractals with infinitely many (countable or even uncountable) similarity maps is not new: cf. \cite{Hata85,Lewe93,Sece12,ArJS17,MaPe17}. They all need an ambient space $\O$. In particular, in \cite{Lewe93} $\O$ is a complete metric space, and a continuous family $\{ f(y,\cdot) : y\in Y\}$ of contractions of $\O$, parametrized by a compact space $Y$, are considered, and a unique fixed-point, which is a self-similar set, is obtained. In \cite{ArJS17} some improvements have been obtained, in case $\O,Y$ are compact metric spaces. Finally, in \cite{Sece12}, if $\O$ is a topological ambient space, $Y=\bn$ and the continuous maps $f_n:\O\to\O$ satisfy some other conditions, a fixed-point is obtained, but the uniqueness can't be assured. 

We conclude by observing that our notion of topological fractal can also be related to the notion of geometric space {\it a la } Whitney: indeed, when $X_0$ is a fixed point, that is to say a topological self-similar fractal, for any $y\in Y$, the map $\pi_y: x_0 \in X_0\mapsto \pi(y,x_0) \in X_1$ may be seen as (the inverse of) a local chart, the Euclidean space  being replaced by the fractal itself, namely the fractal is locally modeled on itself. Transition maps are also well defined, however their domain and range do not have a specific form: they are not homeomorphic to the fractal in general, and they are often finite sets, therefore their regularity doesn't play an important role: in our definition they are simply continuous maps. Nonetheless, overlapping sets, namely $\pi_y(X_0)\cap\pi_{y'}(X_0)$, $y,y'\in Y$, which allow the definition of the transition maps, are crucial, since their knowledge is equivalent to the definition of $\pi$.


\section{Preliminaries on similarity schemes}
We start this section giving the main definition of this paper:
\begin{dfn}[Abelian similarity scheme]
An abelian similarity scheme consists of the following diagram
\begin{equation}\label{SimScheme}
\begin{matrix}
& &Y\times X_0
\\
& &\downarrow\pi\\
X_0&\stackrel{\f}{\hooklongrightarrow}&X_1,
\end{matrix}
\end{equation}
where $X_0,X_1,Y$ are compact Hausdorff spaces, $\f$ is continuous injective, and $\pi$ is continuous surjective. 
\end{dfn}

\begin{prop} An abelian similarity scheme acts on a pair $(Z,\f_Z)$, where $Z$ is compact Hausdorff and $\f_Z:X_0\to Z$ is continuous injective, producing a new pair 
$(\widehat{Z},\f_{\widehat{Z}})=\Phi(Z,\f_{Z})$ such that the following diagram commutes
\begin{equation}\label{SimSchemeAction}
\begin{matrix}
& &Y\times X_0&\stackrel{id_Y \times \f_Z}{\hooklongrightarrow}&Y\times Z
\\
& &\downarrow\pi& &\downarrow\pi_Z\\
X_0&\stackrel{\f}{\hooklongrightarrow}&X_1&\stackrel{\f_{\widehat{Z},X_1}}{\hooklongrightarrow}&{\widehat Z}\,.
\end{matrix}
\end{equation}
\end{prop}
\begin{proof}
We set $\pi_Z$ to be the quotient map  w.r.t. the minimal equivalence relation $\sim$ on $Y\times Z$ such that $\pi(y, x_0)=\pi(y', x'_0)\Rightarrow(y, \f_Z(x_0))\sim (y', \f_Z(x'_0))$. In this way, the map $\f_{Z',X_1}$ is well defined, and we set $\f_{Z'}= \f_{Z',X_1} \circ\f$. 
\end{proof}
\begin{rem}
Let us recall that the push forward of an equivalence relation is not an equivalence relation in general, indeed reflexivity and transitivity may fail. However,
since $\f_Z$ is injective, also $id_Y \times \f_Z$ is injective, therefore transitivity holds. But since $\f_Z$ is not surjective, reflexivity should be imposed. As a consequence, $(y,z)\sim(y',z')$ either if $y=y'$ and $z=z'$ or if there exist $x_0,x_0'$ such that $z=\f_Z(x_0)$, $z'=\f_Z(x_0')$, and $\pi(y,x_0)=\pi(y',x_0')$.
\end{rem}

\begin{prop}
The map $\Phi:(Z,\f)\to({\widehat Z},\f_{\widehat{Z}})$ produced by the action of the similarity scheme is an endo-functor of the category $\ct_{X_0}$ whose objects are pairs $(Z,\f_Z)$, with $Z$ compact Hausdoff and $\f_Z:X_0\to Z$ continuous injective, and whose morphisms $f_\vartheta:(Z,\f_Z)\to(W,\f_W)$ are given by a continuous injective map $\vartheta:Z\to W$ such that $\f_W=\vartheta\circ\f_Z$.
\end{prop}

Isomorphisms in the category are morphisms where $\vartheta$ is an homeomorphism.

\begin{dfn}[Topological self-similar fractal]
Given a similarity scheme $X_0\hookrightarrow X_1\leftarrow Y\times X_0$, we say that a fixed point of the endofunctor $\Phi$ is a topological self-similar fractal w.r.t. the given similarity scheme:
$$
\Phi((Z,\f_Z))\cong(Z,\f_Z).
$$
\end{dfn}

Letting the scheme act on  the pair $(X_1,\f)$ we obtain a new pair $(X_2,\f_{2,0})=\Phi(X_1,\f)$. As mentioned in the introduction, a
repeated application of this procedure gives rise to an infinite diagram
\begin{equation}\label{diagramA}
\hskip-0.8cm \begin{matrix}
&                                                                 
&Y\times X_0&\stackrel{id_Y \times \f_{1,0}}{\hooklongrightarrow}
&Y\times X_1&\stackrel{id_Y \times \f_{2,1}}{\hooklongrightarrow}
&Y\times X_2&\stackrel{id_Y \times \f_{3,2}}{\hooklongrightarrow}
&Y\times X_3&\stackrel{id_Y \times \f_{4,3}}{\hooklongrightarrow}
\\
& 
&\downarrow\pi_0& 
&\downarrow\pi_1& 
&\downarrow\pi_2& 
&\downarrow\pi_3& 
\\
   X_0&\stackrel{\f_{1,0}}{\hooklongrightarrow}
&X_1&\stackrel{\f_{2,1}}{\hooklongrightarrow}
&X_2&\stackrel{\f_{3,2}}{\hooklongrightarrow}
&X_3&\stackrel{\f_{4,3}}{\hooklongrightarrow}
&X_4&\stackrel{\f_{5,4}}{\hooklongrightarrow}
\end{matrix}
\end{equation}
where all the squares are commutative, all vertical maps are surjective, and all horizontal maps  are injective.


\subsection{The extended diagram}
We first observe that by suitably combining the diagrams obtained from \eqref{diagramA} by multiplying by powers of $Y$, we get an infinite stairway, where $\pi_{1,n}:=\pi_n$,
\begin{equation}\label{diagramB}
\hskip-0.8cm \begin{matrix}
& & & & & & &                                                            
&Y^4\times X_0&\stackrel{ id_{Y^4} \times \f_{1,0} }{\hooklongrightarrow}
\\
& & & & & & &
&\downarrow_{ id_{Y^3} \times \pi_{1,0} }& 
\\
& & & & &                                                                
&Y^3\times X_0&\stackrel{ id_{Y^3} \times \f_{1,0} }{\hooklongrightarrow}
&Y^3\times X_1&\stackrel{ id_{Y^3} \times \f_{2,1} }{\hooklongrightarrow}
\\
& & & & &
&\downarrow_{ id_{Y^2} \times \pi_{1,0} }& 
&\downarrow_{ id_{Y^2} \times \pi_{1,1} }& 
\\
& & &                                                                 
&Y^2\times X_0&\stackrel{ id_{Y^2} \times \f_{1,0} }{\hooklongrightarrow}
&Y^2\times X_1&\stackrel{ id_{Y^2} \times \f_{2,1} }{\hooklongrightarrow}
&Y^2\times X_2&\stackrel{ id_{Y^2} \times \f_{3,2} }{\hooklongrightarrow}
\\
& & &
&\downarrow_{ id_{Y} \times \pi_{1,0} }& 
&\downarrow_{ id_{Y} \times \pi_{1,1} }& 
&\downarrow_{ id_{Y} \times \pi_{1,2} }& 
\\
&                                                                 
&Y\times X_0&\stackrel{ id_Y \times \f_{1,0} }{\hooklongrightarrow}
&Y\times X_1&\stackrel{ id_Y \times \f_{2,1} }{\hooklongrightarrow}
&Y\times X_2&\stackrel{ id_Y \times \f_{3,2} }{\hooklongrightarrow}
&Y\times X_3&\stackrel{ id_Y \times \f_{4,3} }{\hooklongrightarrow}
\\
& 
&\downarrow\pi_{1,0}& 
&\downarrow\pi_{1,1}& 
&\downarrow\pi_{1,2}& 
&\downarrow\pi_{1,3}& 
\\
   X_0&\stackrel{\f_{1,0}}{\hooklongrightarrow}
&X_1&\stackrel{\f_{2,1}}{\hooklongrightarrow}
&X_2&\stackrel{\f_{3,2}}{\hooklongrightarrow}
&X_3&\stackrel{\f_{4,3}}{\hooklongrightarrow}
&X_4&\stackrel{\f_{5,4}}{\hooklongrightarrow}
\end{matrix}
\end{equation}
By composing the vertical maps, we get projection maps $\pi_{p,q}:Y^p\times X_q\to X_{p+q}$, defined inductively by $\pi_{1,j}=\pi_j$, $\pi_{p,q}=\pi_{p-1,q+1}\circ(id_{Y^{p-1}} \times \pi_{1,q})$, $p\geq2$.
Let us observe that we also get
\begin{equation}\label{p-q-equation}
\pi_{p,q}=\pi_{p-j,q+j}\circ(id_{Y^{p-j}},\pi_{j,q}), \quad j=1,\dots p-1.
\end{equation}
%


\begin{lem} \label{equivalence} {\bf (Basic equivalence Lemma)} Let $(y_1,\ldots,y_m,\xi_{n-m}) \not= (y'_1,\ldots,y'_m,\xi'_{n-m})$ in $ Y^m\times X_{n-m}$ be such that they provide the same point in $X_n$, i.e.
$$
\pi_{m,n-m}(y_1,\ldots,y_m,\xi_{n-m}) =
\pi_{m,n-m}(y'_1,\ldots,y'_m,\xi'_{n-m})\,.
$$
Then there exists $k\in \{1,\ldots,m\}$, $\xi_0,\xi'_0\in X_0$ such that
\begin{itemize}
\item[$(1)$] $(y_1,\ldots,y_{k-1})=(y'_1,\ldots,y'_{k-1})$ (void if $k=1$)
\item[$(2)$] $\pi_{m-k,n-m}(y_{k+1},\ldots,y_m,\xi_{n-m})=\varphi_{n-k,0}(\xi_0)$ in $X_{n-k}$ ($\xi_{n-m}=\varphi_{n-m,0}(\xi_0)$ if $k=m$)
\item[$(2')$] $\pi_{m-k,n-m}(y'_{k+1},\ldots,y'_m,\xi'_{n-m})=\varphi_{n-k,0}(\xi'_0)$  in $X_{n-k}$ ($\xi'_{n-m}=\varphi_{n-m,0}(\xi'_0)$ if $k=m$)
\item[$(3)$] $(y_k,\xi_0)\sim (y'_k,\xi'_0)$, i.e. $ \pi_{1,0}(y_k,\xi_0)=\pi_{1,0}(y'_k,\xi'_0)$ in $X_0$.
\end{itemize}
\end{lem}
\begin{proof} By induction on $m$.

For $m=1$, by construction of $X_n$ as a quotient of $Y\times X_{n-1}$\,,
$\pi_{1,n-1}(y_1,\xi_{n-1})=\pi_{1,n-1}(y'_1,\xi'_{n-m})$ with $(y_1,\xi_{n-1})\not=(y'_1,\xi'_{n-m})$ means exactly $\xi_{n-1}=\varphi_{n-1}(\xi_0)$ and $\xi'_{n-1}=\varphi_{n-1}(\xi'_0)$ for some $\xi_0,\xi'_0\in X_0$ such that $(y_1,\xi_0)\sim (y'_1,\xi'_0)$.

\smallskip Suppose that the property holds true for $m-1$ and let $(y_1,\cdots,y_m,\xi_{n-m})\not= (y'_1,\cdots,y'_m,\xi'_{n-m})$ in $ Y^m\times X_{n-m}$ have the same projection on $X_n$. Set $\xi_{n-1}=\pi_{m-1,n-m}(y_2,\ldots,y_m,\xi_{n-m})$ and $\xi'_{n-1}=\pi_{m-1,n-m}(y'_2,\ldots,y'_m,\xi'_{n-m})$. There are two possibilities :

Either $(y_1,\xi_{n-1})\not=(y'_1,\xi'_{n-1})$\,: the analysis above of the case $m=1$ provides the result with $k=1$\,.

Or $(y_1,\xi_{n-1})=(y'_1,\xi'_{n-1})$ and the recurrence hypothesis applied to $(y_2,\ldots,y_m,\xi_{n-m})$ and $(y'_2,\ldots,y'_m,\xi'_{n-m})$ provides the result for some $k\in [2,\ldots,m]$.
\end{proof}

\begin{lem}\label{Xk} 
With the notation of Lemma \ref{equivalence} above, one has 
$$
\pi_{m,n-m}(y_1,\ldots,y_m,\xi_{n-m})\in \varphi_{n,k}(X_k).
$$
\end{lem} 
\begin{proof} 
By commutativity of the extended diagram \ref{diagramB}, one has 
\begin{equation*}\begin{split}
\pi_{m,n-m}(y_1,\ldots,y_m,\xi_{n-m})&=\pi_{k,n-k}(y_1,\ldots,y_k,\varphi_{n-k,0}(\xi_{0}))\\
&=\varphi_{n,k}\big(\pi_{k,0}(y_1,\ldots,y_k,\xi_0)\big)\,.
\end{split}\end{equation*}
\end{proof} 

As a Corollary, we get the following result\,:

\begin{lem}\label{neonato} 

Let $x_n\in X_n\backslash \varphi_{n,n-1}(X_{n-1})$ and $(y_1,\ldots,y_n,x_0)$, $(y'_1,\ldots,y'_n,x'_0)$ in $Y^n\times X_0$ both projecting on $x_n$, i.e. such that
$$
x_n=\pi_{n,0}(y_1,\ldots,y_n,x_0)=\pi_{n,0}(y'_1,\ldots,y'_n,x'_0) \in X_n.
$$

\smallskip
\begin{itemize}
\item[$(1)$] Then one has $\displaystyle \left\{\begin{matrix}
(y_1,\ldots,y_{n-1})&=&(y'_1,\ldots,y'_{n-1}) &\text{ in } Y^{n-1} \\
\pi_{1,0}(y_n,x_0)&=&\pi_{1,0}(y'_n,x'_0) &\text{ in } X_1\,.
\end{matrix} \right.$

\item[$(2)$] Stated otherwise, there exists one and only one  $(y_1,\ldots,y_{n-1},x_1)\in Y^{n-1}\times X_1$ whose image by $\pi_{n-1,1}$ is $x_n$.

\item[$(3)$] As a consequence, there exists one and only one pair $(y,\xi_{n-1})\in Y\times X_{n-1}$ such that $\pi_{1,n-1}(y,\xi_{n-1})=x_n$.
\end{itemize} 
\end{lem}
\begin{proof} 
According to Lemmas \ref{equivalence}, \ref{Xk}, if $x_n\not\in \varphi_{n,n-1}(X_{n-1})$, the only possibility in the conclusion of Lemma \ref{equivalence} is $k=n$.
\end{proof}

\begin{lem}\label{X0} Let $1\leq m\leq n$, $x_0\in X_0$, $y_1,\ldots,y_m\in Y$, $\xi_{n-m}\in X_{n-m}$ be such that
$$ 
\pi_{m,n-m}(y_1,\ldots,y_m,\xi_{n-m})=\varphi_{n,0}(x_0)\,.
$$
Then $\xi_{n-m}\in \varphi_{n-m,0}(X_0)$\,.
\end{lem}
\begin{proof} By induction on $n$. 
Case $n=1$\,: then $m=1$ and $\xi_{n-m}\in X_0$.

Suppose now the property holds true for $n-1$ ($n\geq 2$) and $1\leq m \leq n-1$. There exists $(y,\xi_0)\in Y\times X_0$ such that $\varphi_{1,0}(x_0)=\pi_{1,0}(y,\xi_0)$ in $X_1$. We have then 
$$
\varphi_{n,0}(x_0)=\varphi_{n,1} \circ \varphi_{1,0}(x_0)=\varphi_{n,1}(\pi_{1,0}(y,\xi_0))=\pi_{1,n-1}(y,\varphi_{n-1,0}(\xi_0))\,.
$$
Setting $\xi_{n-1}=\pi_{m-1,n-m}(y_2,\ldots,y_m,\xi_{m-n})$, we see that $(y,\varphi_{n-1}(\xi_0))$ and $(y_1,\xi_{n-1})$ in $Y\times X_{n-1}$ have the same projection on $X_n$. There are thus two possibilities\,:

either they are equal, and then $\xi_{n-1}=\varphi_{n-1}(\xi_0)\in \varphi_{n,0}(X_{n-1})$ and the induction hypothesis provides the result;

or they are not equal, and this implies again $\xi_{n-1}\in \varphi_{n-1,0}(X_0)$, so that the induction hypothesis provides the result.
\end{proof}

\begin{lem}\label{Xm} 
Let $1\leq m\leq n$, $x_m\in X_m$ and $(y_1,\ldots, y_m,\xi_{n-m})\in Y^m\times X_{n-m}$ be such that 
$$
\pi_{m,n-m}(y_1,\ldots, y_m,\xi_{n-m})=\varphi_{n,m}(x_m)\,.
$$
Then $\xi_{n-m}\in \varphi_{n-m,0}(X_0)$\,.
\end{lem}
\begin{proof} There exists $(y'_1,\ldots,y'_m,\xi'_0)\in Y^m\times X_0$ such that $\pi_{m,0}(y'_1,\ldots,y'_m,\xi'_0)=x_m$ in $X_m$. Hence $(y_1,\ldots, y_m,\xi_{n-m})$ and $(y'_1,\ldots,y'_m,\varphi_{n-m,0}(\xi'_0))$ have the same projection on $X_n$.

Case 1. $(y_1,\ldots, y_m,\xi_{n-m})=(y'_1,\ldots,y'_m,\varphi_{n-m,0}(\xi'_0))$\,: the conclusion holds true.

Case 2. $(y_1,\ldots, y_m,\xi_{n-m})\not=(y'_1,\ldots,y'_m,\varphi_{n-m,0}(\xi'_0))$\,: by Lemma \ref{equivalence},  there exist $k\in [1,\ldots,m]$ such that $\pi_{m-k,n-m}(y_{k+1},\ldots,y_m, \xi_{n-m})\in \varphi_{n-k,0}(X_0)$. By Lemma \ref{X0} the conclusion holds true.
\end{proof}

\begin{cor}\label{Xn+1} 
If $(y,\xi_n)\in Y\times X_n$ is such that $\pi_{1,n}(y,\xi_n)\in \varphi_{n+1,n}(X_n)$, then $\xi_n\in \varphi_{n,n-1}(X_{n-1})$. 
\end{cor}
\begin{proof} 
Write $\xi_n=\pi_{n,0}(y_1,\ldots,y_n,x_0)$ for some suitable $(y_1,\ldots,y_n,x_0)\in Y^n\times X_0$. As $\pi_{n+1,0}(y,y_1,\ldots,y_n,x_0)=\pi_{1,n}(y,\xi_n)$ lies in $\varphi_{n+1,n}(X_n)$, Lemma \ref{Xm} provides $\pi_{1,0}(y_n,x_0)=\varphi_{1,0}(\xi_0)$ for some $\xi_0\in X_0$, hence 
$$
\xi_n=\pi_{n-1,1}(y_1,\ldots,y_{n-1},\varphi_{1,0}(\xi_0))=\varphi_{n,n-1}\big(\pi_{n-1,0}(y_1,\ldots,y_{n-1},\xi_0)\big)\in \varphi_{n,n-1}(X_{n-1}).
$$
\end{proof}

\section{A combinatorial equivalence relation}\label{section:CombDefXinfinity}

By using the map $\pi_{n,0}$, any point $\vec{y}_n\in Y^n$ gives rise to the set \begin{equation}\label{cell}
C(\vec{y}_n):=\pi_{n,0}(\{(\vec{y}_n,x):x\in X_0\})\subset X_n.
\end{equation}
Such set is called the {\it cell} determined by $\vec{y}_n$. For convenience, we set $Y^0:=\{\emptyset\}$, $\pi_{0,0}:=id_{X_0}$, and $C(\emptyset):=\{\pi_{0,0}(x):x\in X_0\}=X_0$.

\begin{dfn}
Given $x_n\in X_n$, we define $\Gamma_n(x_n)\subset Y^\bn$ as
$$
\Gamma_n(x_n)=\{y_\infty\in Y^\bn:\f_{p,n}(x_n)\in C(\rho_p(y_\infty)),\forall p\geq n\},
$$
where $\rho_p:Y^\bn\to Y^p$ denotes the $p$-th truncation map.
\end{dfn}

Observe that $\Gamma_n(x_n)$ is obviously a closed subset of $Y^\N$.
\\
By a diagram chasing argument on the diagram \eqref{diagramB}, it is easy to see that $\Gamma_n(x_n)$ is not empty.
\\

\begin{prop}\label{proj1} 
For $n\in \N$ and $x_n\in X_n$, we have $\Gamma_{n+1}(\varphi_{n+1,n}(x_n))=\Gamma_n(x_n)$.
\end{prop}
\begin{proof} 
We have $y_\infty\in \Gamma_{n+1}(\varphi_{n+1,n}(x_n))$ iff $\varphi_{n+k,n}(x_n)\in C(y_1,\ldots,y_{n+k})$ for all $k\geq 1$. All we have to prove is that $y_\infty\in \Gamma_{n+1}(\varphi_{n+1,n}(x_n))$ implies $x_n\in C(y_1,\ldots,y_n)$. 

If $y_\infty\in \Gamma_{n+1}(\varphi_{n+1,n}(x_n))$, then $\varphi_{n+1,n}(x_n)\in C(y_1,\ldots,y_{n+1})$, i.e. $\exists \xi_0\in X_0$ such that $\varphi_{n+1,n}(x_n)=\pi_{n+1,0}(y_1,\ldots,y_{n+1},\xi_0)$.

Setting $\xi_1=\pi_{1,0}(y_{n+1},\xi_0)\in X_1$, we have $\varphi_{n+1,n}(x_n)=\pi_{n,1}(y_1,\ldots,y_{n},\xi_1)$. 

Hence $\pi_{n+1,0}(y_1,\ldots,y_{n},\xi_1)\in \varphi_{n+1,n}(X_n)$ which, by Lemma \ref{Xm}, implies $\xi_1\in \varphi_{1,0}(X_0)$, i.e. $\xi_1=\varphi_{1,0}(\xi'_0)$ for some $\xi'_0\in X_0$.

This means
\begin{equation*}\begin{split}
\varphi_{n+1,n}(x_n)&=\pi_{n,1}\big(y_1,\ldots,y_n,\varphi_{1,0}(\xi'_0)\big) \\
&=\varphi_{n+1,n}\big(\pi_{n,0}(y_1,\ldots,y_n,\xi'_0)\big)
\end{split}\end{equation*}
and, by injectivity of $\varphi_{n+1,n}$, $x_n=\pi_{n,0}(y_1,\ldots,y_n,\xi'_0)$, which ends the proof.
\end{proof}

\begin{lem}\label{extension} 
Let $x_n=\pi_{n,0}(y_1,\ldots,y_n,x_0)\in X_n$, and $y_\infty'\in \Gamma_0(x_0)$. Then the element $(y_1,\ldots,y_n,y'_1,\ldots, y'_m,\ldots)$ of $Y^\N$ lies in $\Gamma_n(x_n)$.

In particular, this shows that for any $(y_1,\ldots,y_n)\in Y^n$, there exists at least one $y_\infty\in \Gamma_n(x_n)$, of which the  first $n$ coordinates are exactly $y_1,\ldots,y_n$.
\end{lem}
\begin{proof} 
Straightforward. 
\end{proof}

Notice that this Lemma proves also $\Gamma_n(x_n)\not=\emptyset$, since for $x_0\in X_0$, we can choose $(y_1,\xi_0^1)\in Y\times X_0$ such that $\pi_{1,0}(y_1,\xi_0^1)=\varphi_{1,0}(x_0)$, then $(y_2,\xi_0^2)\in Y\times X_0$ such that $\pi_{1,0}(y_2,\xi_0^2)=\varphi_{1,0}(\xi_0^1)$ hence $\pi_{2,0}(y_1,y_2,\xi_0^2)=\varphi_{2,0}(x_0)$, and so on, in order to get an element of $\Gamma_0(x_0)$.

\begin{lem}\label{Gammaneonato} 
Let $n\geq 2$, $x_n\in X_n\backslash \varphi_{n,n-1}(X_{n-1})$, and $(y_1,\ldots,y_{n-1},x_1)$ be the only element of $Y^{n-1}\times X_1$ which is projected on $x_n$ by $\pi_{n,1}$. Then one has 

\begin{itemize}
\item[$(1)$] $\Gamma_n(x_n) \supset \{ \big(y_1,\ldots,y_{n-1}, y_\infty) : y_\infty \in \Gamma_1(x_1) \}$,

\item[$(2)$] if $\pi_{1,0}(y,\cdot) : X_0 \to X_1$ is injective, $\Gamma_n(x_n) = \{ \big(y_1,\ldots,y_{n-1},y_\infty) : y_\infty \in \Gamma_1(x_1) \}$.
\end{itemize}
\end{lem}
\begin{proof} 
Apply Lemma \ref{neonato}.
\end{proof}

Given two elements $y_\infty,y'_\infty\in Y^\bn$, consider the relation $\car$ on $Y^\bn$ defined as follows:
\begin{equation}\label{car}
y_\infty\, \car\, y'_\infty \mathrm{\, if\, } 
\begin{cases}
\mathrm{\, either\, }& y'_\infty=y_\infty,\\
\mathrm{\, or\, } & y_\infty\ne y'_\infty\, \&\,\exists n\in\bn\cup\{0\}, x_n\in X_n: y_\infty,y'_\infty\in\Gamma(x_n).
\end{cases}
\end{equation}

We stress that $\car$ is not an equivalence relation in general, since transitivity is not necessarily satisfied. 

\begin{dfn}
Let us  define $\car_\infty$ as the minimal equivalence relation containing $\car$, set $X_\infty=Y^\bn/\car_\infty$, and denote by $\pi^\bn:Y^\bn\to X_\infty$ the associated projection map.
\end{dfn}

Let us remark that the pair $(n,x_n)$ in \eqref{car} may not be unique. 

We are now able to prove Property \ref{embeddings} mentioned in the Introduction.
\begin{thm}\label{embeddings2}
There exist maps $\f_{\infty,n}:X_n\to X_\infty$ such that
\begin{itemize}
\item[a)] $\f_{\infty,n}\circ\f_{n,m}=\f_{\infty,m}$
\item[b)] $\bigcup_n \varphi_{\infty,n}(X_n)$ is dense in $X_\infty$.
\end{itemize}
\end{thm}
\begin{proof} 
(a) Since two elements of $\Gamma(x_n)$ are related by $\car$, a fortiori $\pi^\bn(\Gamma(x_n))$ is a singleton in $X_\infty$, whose unique element will be denoted by $\f_{\infty,n}(x_n)$. 
The maps  are clearly compatible with the maps $\f_{n,m}:X_m\to X_n$ obtained by combining the maps $\f_{p+1,p}$ of diagram \eqref{diagramA}.
\\
(b) Straightforward by Lemma \ref{extension}.
\end{proof}

\subsection{Fully injective similarity schemes}

\begin{dfn}
We say that a similarity scheme is {\it transitive} if the relation $\car$ satisfies the transitive property, and that a scheme is {\it fully injective} if the maps $\f_{\infty,n}$ are injective.
\end{dfn}

Let us consider the following three properties:
\begin{enumerate}
\item[(P1)] The maps $\f_{\infty,n}$ are injective (full injectivity);
\item[(P2)] If for some $n\in\bn$, $y_\infty\in Y^\bn$, $\exists x_n\in X_n$ such that $y_\infty\in\Gamma(x_n)$, such $x_n$ is unique.
\item[(P3)] If for some $n\geq1$, $x_n\in X_n\setminus\f_{n,n-1}(X_{n-1})$, $y_\infty\in Y^\bn$, we have $y_\infty\in\Gamma(x_n)$, the pair $(n,x_n)$ is unique.
\end{enumerate}%

\begin{prop} \label{trans}
Let $X_0 \rightarrow X_1 \leftarrow Y \times X_0$ be a similarity scheme. If property $(P2)$ is satisfyed then the scheme is transitive and fully injective (i.e. {\rm (P1)} holds true).
\end{prop}
\begin{proof}
By the uniqueness assumption it is easy to see that $\car$ is an equivalence relation, indeed if $y_\infty\car y'_\infty$ and $y'_\infty\car y''_\infty$ there exist $x_n\in X_n$, $x_m\in X_m$ such that $y_\infty,y'_\infty\in\Gamma_n(x_n)$ and $y'_\infty,y''_\infty\in\Gamma_m(x_m)$. It is not restrictive to assume that $m\geq n$, therefore $y_\infty,y'_\infty\in\Gamma_m(\f_{m,n}(x_n))$. Uniqueness implies that $x_m=\f_{m,n}(x_n)$, hence $y_\infty,y'_\infty,y''_\infty\in\Gamma_m(x_m)$.
\\
Assume now $x_n,x'_n\in X_n$ verify $\f_{\infty,n}(x_n)=\f_{\infty,n}(x'_n)$. Then there exists $y_\infty,y'_\infty\in Y^{\bn}$ such that $y_\infty\in\Gamma_n(x_n)$, $y'_\infty\in\Gamma_n(x'_n)$, and $y_\infty\,\car\,y'_\infty$, where we used the transitive property of $\car$. Then, either $y_\infty=y'_\infty$, in which case $x_n=x'_n$ by the uniqueness hypothesis, or, by \eqref{car}, $\exists m\in\bn, x_m\in X_m$ such that $y_\infty,y'_\infty\in \Gamma_m(x_m)$. Reasoning as above, e.g. for $m\geq n$, $\f_{m,n}(x_n)=\f_{m,n}(x'_n)=x_m$, which implies $x_n=x'_n$ by the injectivity of $\f_{m,n}$. 
\end{proof}


\begin{prop}\label{fullness}
Let $X_0 \rightarrow X_1 \leftarrow Y \times X_0$ be a similarity scheme. Then
properties $(P1)$, $(P2)$ and $(P3)$ are equivalent:
\end{prop}
\begin{proof}
$(P1)\Rightarrow(P2)$. Assume that for some $n\in\bn$, $y_\infty\in Y^\bn$, $x_n,x'_n\in X_n$ we have $y_\infty\in\Gamma_n(x_n)\cap\Gamma(x'_n)$, which implies that $\pi^\bn(\Gamma_n(x_n))=\pi^\bn(\Gamma_n(x'_n))$, namely $\f_{\infty,n}(x_n)=\f_{\infty,n}(x'_n)$. By the injectivity hypothesis, $x_n=x'_n$.
\\
$(P2)\Rightarrow(P3)$.
Assume we have $(n,x_n), (m,x_m)$, $m\geq n$, $y_\infty\in Y^\bn$ such that $y_\infty\in \Gamma_n(x_n)\cap\Gamma_m(x_m)$. By Proposition \ref{proj1}, $y_\infty\in\Gamma(\f_{m,n}(x_n))$, hence, by hypothesis, $x_m=\f_{m,n}(x_n)$. If $x_n\in X_n\setminus\f_{n,n-1}(X_{n-1})$, then $(n,x_n)$ is the unique pair. If not, $\exists! p<n$ such that $x_p\in X_p\setminus\f_{p,p-1}(X_{p-1})$ with $\f_{n,p}(x_p)=x_n$. By Proposition \ref{proj1}, $\Gamma_n(x_n)=\Gamma_n(\f_{n,p}(x_p))=\Gamma_p(x_p)$, hence $y_\infty\in\Gamma_p(x_p)$, so $(p,x_p)$ is the unique pair. 
\\
$(P3)\Rightarrow(P2)$. Assume that, for  given $n\in\bn$, we get $x_n\in X_n$ and $y_\infty\in Y^\bn$ such that $y_\infty\in\Gamma_n(x_n)$. Clearly we find a unique $p\leq n$, and $x_p\in X_p\setminus X_{p-1}$ such that $x_n=\f_{n,p}(x_p)$. By Proposition \ref{proj1}, $\Gamma_n(x_n)=\Gamma_n(\f_{n,p}(x_p))=\Gamma_p(x_p)$. Since $x_p$ is unique by hypothesis, uniqueness for $x_n$ follows.
\\
$(2)\Rightarrow(1)$. See Proposition \ref{trans}.
\end{proof}

%
%
%
%
%

\subsection{Discrete similarity schemes} 

\begin{dfn}\label{discrete} 
The similarity scheme $X_0 \rightarrow X_1 \leftarrow Y \times X_0$ is said to be {\it discrete} if each cell in $X_1$ intersects $\varphi_{1,0}(X_0)$ in at most one point.
\end{dfn}
 
Otherwise stated, for fixed $y\in Y$, if $x_0,x'_0\in X_0$ are such that $\pi_{1,0}(y,x_0), \pi_{1,0}(y,x'_0)\in \varphi(X_0)$, then  $\pi_{1,0}(y,x_0)=\pi_{1,0}(y,x'_0)$ in $X_1$.

\begin{lem} \label{Hn} 

If $X_0 \rightarrow X_1 \leftarrow Y \times X_0$ is discrete, then for every $n\geq 1$, each cell $C(\vec{y}_n) \subset X_n$ intersects $\varphi_{n,n-1}(X_{n-1})$ in at most one point.
\end{lem}
Otherwise stated, for fixed $y_1,\ldots y_n \in Y_n$, if $x_0,x'_0 \in X_0$ are such that $\pi_{n,0}(y_1,\ldots,y_n,x_0)$ and $\pi_{n,0}(y_1,\ldots,y_n,x'_0) \in \varphi_{n,n-1}(X_{n-1})$, then $\pi_{n,0}(y_1,\ldots,y_n,x_0)=\pi_{n,0}(y_1,\ldots,y_n,x'_0)$.
\begin{proof} 
By induction on $n$. Suppose we have already proved the statement for some $n\geq 2$, and let $y_1,\ldots,y_n\in Y$, $x_0,x'_0\in X_0$, and $\xi_{n-1},\xi'_{n-1}\in X_{n-1}$ be such that $\pi_{n,0}(y_1,\ldots,y_n,x_0)= \varphi_{n,n-1}(\xi_{n-1})$ and $\pi_{n,0}(y_1,\ldots,y_n,x'_0) = \varphi_{n,n-1}(\xi'_{n-1})$ in $\varphi_{n,n-1}(X_{n-1})$.

Then $\pi_{1,n-1}(y_1,\pi_{n-1,0}(y_2,\ldots,y_n,x_0)) = \pi_{n,0}(y_1,\ldots,y_n,x_0)= \varphi_{n,n-1}(\xi_{n-1})$, so, by Corollary \ref{Xn+1}, $\pi_{n-1,0}(y_2,\ldots,y_n,x_0) \in \f_{n-1,n-2}(X_{n-2}) \cap C(y_2,\ldots,y_n)$, and therefore it is unique, by the inductive hypothesis. Analogously $\pi_{n-1,0}(y_2,\ldots,y_n,x_0')$ is the unique element of $\f_{n-1,n-2}(X_{n-2}) \cap C(y_2,\ldots,y_n)$, so it coincides with $\pi_{n-1,0}(y_2,\ldots,y_n,x_0)$.
 
%
%
%

Finally, 
\begin{equation*}\begin{split}
\pi_{n,0}(y_1,\ldots,y_n,x'_0)&=\pi_{1,n-1}\big(y_1,\pi_{n-1,0}(y_2,\ldots,y_n,x'_0)\big)\\
&=\pi_{1,n-1}\big(y_1,\pi_{n-1,0}(y_2,\ldots,y_n,x_0)\big) \\
&=\pi_{n,0}(y_1,\ldots,y_n,x_0),
\end{split}\end{equation*}
which ends the proof. 
\end{proof}

\begin{prop}\label{discrete2} 
If the similarity scheme is discrete, then it is fully injective, hence transitive.
\end{prop}
\begin{proof} 
Suppose that there exist $n\geq 0$, $x_n,x'_n\in X_n$ and $y_\infty\in Y^\N$ such that $\forall k\geq 0$, $\varphi_{n+k,n}(x_n),\varphi_{n+k,n}(x'_n)\in C(y_1,\ldots,y_{n+k})\,.$
  
Write this property for $k=1$\,: $\varphi_{n+1,n}(x_n)$ and $\varphi_{n+1,n}(x'_n)$ are two points of $C(y_1,\ldots, y_n,y_{n+1})$ belonging to $X_{n+1}\backslash \varphi_{n+1,n}(X_n)$. By Lemma \ref{Hn} they must be equal and, as $\varphi_{n+1,n}$ is injective, we get $x'_n=x_n$.
\end{proof}




\section{Properties of the quotient space $X_\infty$.}

From now on we shall assume that the similarity scheme is fully injective.

\subsection{The subsets $\Gamma(x_n)$.}

\begin{prop}\label{cond} 
Let $n,m\in \N$,  $x_n\in X_n$ and $x'_m\in X_m$.

\item[(1)] The following conditions are equivalent\,:

\begin{itemize}

\item[(i)] $\left\{ \begin{matrix} \text{ either } & n=m & \text{ and } & x_n=x'_m \\ \text{ or } &m<n & \text{ and } &x_n=\varphi_{n,m}(x'_m) \\ \text{ or } &m>n & \text{ and } &x'_m=\varphi_{m,n}(x_m), 
 \end{matrix} \right.$
 
\item[(ii)] for all $p\geq \max(m,n)$, $\varphi_{p,m}(x'_m)=\varphi_{p,n}(x_n)$,

\item[(iii)] there exists $p\geq \max(m,n)$, such that $\varphi_{p,m}(x'_m)=\varphi_{p,n}(x_n)$,

\item[(iv)] $\varphi_{n+m,n}(x_n)=\varphi_{m+n,m}(x'_m)$\,.
\end{itemize}

\item[(2)] We have $\Gamma_n(x_n)=\Gamma_m(x'_m)$ whenever one of the conditions above is satisfied.
\end{prop}
\begin{proof} 
$(1)$ is obvious, while $(2)$ is a straightforward corollary of Proposition \ref{proj1}. 
\end{proof}

  
\begin{prop}\label{GammainterGamma} 
For $n,m\geq 0$, $x_n\in X_n$, $x'_m\in X_m$, the following alternative holds true :
\begin{itemize} 

\item[-] either $\varphi_{m+n,n}(x_n)=\varphi_{n+m,m}(x'_m)$ and $\Gamma_n(x_n)=\Gamma_m(x'_m)$
\item[-] or $\Gamma_n(x_n)\bigcap \Gamma_m(x'_m)=\emptyset.$
\end{itemize}
  
Otherwise stated, one has either $\Gamma_n(x_n)=\Gamma_m(x'_m)$ or $\Gamma_n(x_n)\bigcap \Gamma_m(x'_m)=\emptyset$, and the equality occurs only in the equivalent situations of Proposition \ref{cond} $(1)$.
\end{prop}
\begin{proof} 
From Proposition \ref{cond}, 
it is enough to show that $\Gamma_n(x_n)\bigcap \Gamma_m(x'_m)\not=\emptyset$ implies $\varphi_{m+n,n}(x_n)=\varphi_{n+m,m}(x'_m)$.

Set $\xi_{n+m}=\varphi_{m+n,n}(x_n)$ and $\xi'_{n+m}=\varphi_{n+m,m}(x'_m)$ in $X_{n+m}$. Suppose $y_\infty\in \Gamma_n(x_n)\bigcap \Gamma_m(x'_m)=\Gamma_{n+m}(\xi_{n+m})\bigcap \Gamma_{n+m}(\xi'_{n+m})$. From Proposition \ref{fullness}, we obtain $\xi_{n+m}=\xi'_{n+m}$.
\end{proof}

Next Lemma states a kind of continuity for the map $x_n\to \Gamma(x_n)$.

\begin{lem}\label{cont1} 
Fix $n\geq 0$. Let $(x_n^{(\a)})_{\a\in\ca}$ be a  net in $X_n$, converging to $x_n\in X_n$, and $(y_\infty^{(\a)})_{\a\in\ca}$ a net in $Y^\N$, converging to $y_\infty\in Y^\bn$. Suppose moreover $y_\infty^{(\a)}\in \Gamma_n(x_n^{(\a)})$, for any $\a\in\ca$. Then $y \in \Gamma_n(x_n)$.
\end{lem}
\begin{proof} 
%
For any $p\geq n$, $\a\in\ca$, there is $\xi_{p,\a}\in X_0$ such that $\f_{p,n}(x_n^{(\a)})=\pi_{p,0}(y_1^{(\a)},\ldots,y_p^{(\a)},\xi_{p,\a})$. Since $X_0$ is compact Hausdorff, upon passing to a subnet, we can assume that $\xi_{p,\a} \to \xi_p\in X_0$. From continuity of $\f_{p,n}$ and $\pi_{p,0}$, we get $\f_{p,n}(x_n)=\pi_{p,0}(y_1,\ldots,y_p,\xi_p)$, for all $p\geq n$, so $y\in \G_n(x_n)$.
\end{proof}

%
%

\begin{thm} \label{closedRel}
$\mathcal R_\infty$ is a closed equivalence relation on $Y^\N$, namely the space $X_\infty$ endowed with the quotient topology is compact Hausdorff.
\end{thm}  
\begin{proof} 
$\mathcal R_\infty$ is an equivalence relation by Proposition \ref{fullness}.

In order to prove that its graph is closed, take $(\overset{\rightarrow}{y}^{(\a)})_{\a\in\ca}$ and $(\overset{\rightarrow}{z}^{(\a)})_{\a\in\ca}$ two convergent nets in $Y^\N$, $\lim_{\a\in\ca} \overset{\rightarrow}{y}^{(\a)}=\overset{\rightarrow}{y}^{(\infty)}$, $\lim_{\a\in\ca} \overset{\rightarrow}{z}^{(\a)}=\overset{\rightarrow}{z}^{(\infty)}$, such that $\overset{\rightarrow}{y}^{(\a)}\,\mathcal R_\infty \overset{\rightarrow}{z}^{(\a)}$ for all $\alpha$, and prove that $\overset{\rightarrow}{y}^{(\infty)} \mathcal R_\infty \overset{\rightarrow}{z}^{(\infty)}$.

Suppose first that the set $\ca_0$ of $\a$'s such that $\overset{\rightarrow}{y}^{(\a)}=\overset{\rightarrow}{z}^{(\a)}$ is cofinal. Then we have $\overset{\rightarrow}{y}^{(\infty)}=\overset{\rightarrow}{z}^{(\infty)}$ and the result follows.

We can thus suppose, without loss of generality, that, for any $\a$, there exists $x_{n(\a)}\in X_{n(\a)}$ such that $\overset{\rightarrow}{y}^{(\a)}$, $\overset{\rightarrow}{z}^{(\a)}\in \Gamma_{n(\a)}(x_{n(\a)})$. Moreover, by using Proposition \ref{proj1}, we can choose $n(\a)$ as small as possible, i.e. such that $x_{n(\a)}\in X_{n(\a)}\backslash X_{n(\a)-1}$. By Lemma \ref{neonato} we have then 
\begin{equation}\label{eqnk}
(y^{(\a)}_1,\ldots, y^{(\a)}_{n(\a)-1})=(z^{(\a)}_1,\ldots, z^{(\a)}_{n(\a)-1}).
\end{equation}

Restricting oneself to a subsequence if necessary, we have only two cases to consider :

-- either $\{n(\a):\a\in\ca\}$ is unbounded, and, upon passing to a subnet, we can assume $\lim_{\a\in\ca} n(\a)=\infty$, in which case  \eqref{eqnk} provides in the limit $\overset{\rightarrow}{y}^{(\infty)} = \overset{\rightarrow}{z}^{(\infty)}$;

-- or $N=\sup_{\a\in\ca} n(\a)$ is finite. In this case, we replace $x_{n(\a)}$ with $\xi^{(\a)}=\varphi_{N,n(\a)}(x_{n(\a)})\in X_N$. By compactness of $X_N$, we can suppose that the net $\xi^{(\a)}$ has a limit $\xi^{(\infty)}$ and, in the limit, by Lemma \ref{cont1} we have $\overset{\rightarrow}{y}^{(\infty)}$, $\overset{\rightarrow}{z}^{(\infty)}\in \Gamma_N(\xi^{(\infty)})$.
\end{proof}

\section{$X_\infty$ as a self-similar space.}
In this Section we prove that the limit space $X_\infty$ is a topologial self-similar fractal, i.e. it is a fixed point of the endo-functor associated with the given similarity scheme.

According to diagram \eqref{FixedPoint}, the map $id_Y\times\f_{\infty,0}:Y\times X_0\to Y\times X_\infty$ gives rise to the diagram
\begin{equation}\label{FixedPointBis}
\begin{matrix}
& &Y\times X_0&\stackrel{ id_Y \times \f_{\infty,0} }{\hooklongrightarrow}&Y\times X_\infty
\\
& &\downarrow\pi& &\downarrow\pi_{1,\infty}\\
X_0&\stackrel{\f}{\hooklongrightarrow}&X_1&\stackrel{\f_{\widehat X_\infty,X_1}}{\hooklongrightarrow}& \widehat X_\infty\,.
\end{matrix}
\end{equation}

Let us recall that, by formula \eqref{SimSchemeAction}, the space $\widehat X_\infty$ is the quotient of $Y\times X_\infty$ by the equivalence relation where $(y,x_\infty)$ and $(y', x'_\infty)$ are equivalent if either they are equal or if there exist $x_0,x'_0\in X_0$ such that $x_\infty=\f_{\infty,0}(x_0)$, $x'_\infty=\f_{\infty,0}(x'_0)$ and  $\pi_{1,0}(y,x_0)=\pi_{1,0}(y',x'_0)$ in $X_1$.

An alternative way is to consider that $\widehat X_\infty$ is the quotient space of $Y\times Y^\N$ by the equivalence relation where $(y, y_\infty)$ is equivalent to $(y',y_\infty')$ if 

-- either $y=y'$ and there exist $n\geq 0$, $x_n\in X_n$ such that $y_\infty,y_\infty'\in \Gamma_n(x_n)$

-- or there exist $x_0,x'_0\in X_0$ such that $y_\infty\in \Gamma_0(x_0)$, $y_\infty'\in \Gamma_0(x'_0)$ and $\pi_{1,0}(y,x_0) = \pi_{1,0}(y',x'_0)$.

A third equivalent description of $\widehat X_\infty$ is to consider it as the quotient of $Y^\N$ by the equivalence relation $\widehat {\mathcal R}_\infty$ thus defined\,:

\centerline{$y_\infty \widehat {\mathcal R}_\infty y_\infty'$ if $\left\{  \begin{matrix} \text{ either }  y_1=y'_1  \text{ and } \exists\, n\geq 0\,,\,x_n\in X_n \text{ s.t. } &\left\{ \begin{matrix} 
 (y_2,\ldots,y_k,\ldots)\in \Gamma_n(x_n)\,, \\ (y'_2,\ldots,y'_k,\ldots)\in \Gamma_n(x_n)\,;\end{matrix}\right. \\
 \text{ or } \exists\,x_0,x'_0\in X_0 \text{ s.t. } \pi_{1,0}(y,x_0) = \pi_{1,0}(y',x'_0) \text{ and } &\left\{ \begin{matrix} 
 (y_2,\ldots,y_k,\ldots)\in \Gamma_0(x_0)\,, \\ (y'_2,\ldots,y'_k,\ldots)\in \Gamma_0(x'_0)\,.\end{matrix} \right.
\end{matrix}\right.$}

\begin{thm}\label{LimFixPt}
One has $\mathcal R_\infty= \widehat {\mathcal R}_\infty$, which means that
$X_\infty$, with the imbedding $\varphi_{\infty,0}:X_0\to X_\infty$, is a fixed point of the similarity scheme. It will be called the limit fixed point.
\end{thm}
\begin{proof} 
We suppose first that $y_\infty\, \widehat {\mathcal R}_\infty y_\infty'$ and we prove $y_\infty\,  {\mathcal R}_\infty y_\infty'$.
\smallskip

Case 1\,: $y_1=y'_1  \text{ and } \exists\, n\geq 0\,,\,x_n\in X_n \text{ s.t. } \left\{ \begin{matrix} 
 (y_2,\ldots,y_k,\ldots)\in \Gamma_n(x_n)\,, \\ (y'_2,\ldots,y'_k,\ldots)\in \Gamma_n(x_n)\,.\end{matrix}\right.$
 
Set $x_{n+1}=\pi_{1,n}(y_1,x_n)$. For any $k\geq 2$, one can find $x_{0,k},x'_{0,k}\in X_0$ such that $x_n=\pi_{n+k-1,0}(y_2,\ldots,y_{n+k},x_{0,k})$, which implies $x_{n+1}=\pi_{n+k,0}(y_1,\ldots,y_{n+k},x_{0,k})$\,: we have shown $y_\infty\in \Gamma_{n+1}(x_{n+1})$. Similarly $y_\infty'\in \Gamma_{n+1}(x_{n+1})$, hence $y_\infty\, \mathcal R y_\infty'$.
 
\smallskip Case 2\,: $\exists\,x_0,x'_0\in X_0 \text{ s.t. } \pi_{1,0}(y_1,x_0) = \pi_{1,0}(y'_1,x'_0) \text{ and } \left\{ \begin{matrix} 
 (y_2,\ldots,y_k,\ldots)\in \Gamma_0(x_0)\,, \\ 
 (y'_2,\ldots,y'_k,\ldots)\in \Gamma_0(x'_0)\,.
\end{matrix} \right.$
Set $x_1=\pi_{1,0}(y_1,x_0)=\pi_{1,0}(y'_1,x'_0)$.
 
For any $k\geq 2$ one can find $\xi_{0,k}\in X_0$ such that $\pi_{k-1,0}(y_2,\ldots,y_k,\xi_{0,k})=\varphi_{k-1,0}(x_0)$, which implies $\pi_{k,0}(y_1,y_2,\ldots,y_k,\xi_{0,k})=\varphi_{k,1}\big(\pi_{1,0}(y_1,x_0)\big)= \varphi_{k,1}(x_1)$. We have proved $y_\infty\in \Gamma_1(x_1)$. Similarly $y_\infty'\in \Gamma_1(x_1)$, hence $y_\infty\,  {\mathcal R}_\infty y_\infty'$.
 
\medskip 

We suppose now that $y_\infty\, {\mathcal R}_\infty y_\infty'$ and we prove $y_\infty\,  \widehat {\mathcal R}_\infty y_\infty'$.
 
 If $y_\infty=y_\infty'$, there is nothing to prove. We can thus suppose $y_\infty,y_\infty'\in \Gamma_n(x_n)$ for some $n\geq 0$ and some $x_n\in X_n$. We shall distinguish two cases.
 
\smallskip
Case 1\,: $n\geq 2$ and $x_n\in X_n\backslash \varphi_{n,n-1}(X_{n-1})$. By Lemma \ref{neonato}, there exists a unique pair $(y,\xi_{n-1})\in Y\times X_{n-1}$ such that $\pi_{1,n-1}(y,\xi_{n-1})=x_n$.  For $k\geq 0$, there exists $\xi_{0,k}\in X_0$ such that $\pi_{n+k,0}(y_1,\ldots, y_{n+k},\xi_{0,k})=\varphi_{n+k,n}(x_n)$, hence $\pi_{1,n+k-1}\big(y_1,\pi_{n+k-1,0}(y_2,\ldots,y_{n+k},\xi_{0,k})\big)=\varphi_{n+k,n}(x_n)$. On the other hand, $\varphi_{n+k,n}(x_n)=\pi_{1,n+k-1}\big(y,\varphi_{n+k-1,n-1}(\xi_{n-1})\big)$. Assuming $\big(y_1,\pi_{n+k-1,0}(y_2,\ldots,y_{n+k},\xi_{0,k})\big) \not=\big(y,\varphi_{n+k-1,n-1}(\xi_{n-1})\big)$ would imply $\varphi_{n+k-1,n-1}(\xi_{n-1})\in \varphi_{n+k-1,0}(X_0)$ and $x_n\in \f_{n,1}(X_1) \subset \f_{n,n-1}(X_{n-1})$, which would lead to a contradiction.

We have shown $\big(y_1,\pi_{n+k-1,0}(y_2,\ldots,y_{n+k},\xi_{0,k})\big)=\big(y,\varphi_{n+k-1,n-1}(\xi_{n-1})\big)$, from which one deduces $y_1=y$ and $\varphi_{n+k-1,n-1}(\xi_{n-1})\in C(y_2,\ldots,y_{n+k})$ for all $k$. Finally, we have $y_1=y$ and $(y_2,\ldots,y_k, \ldots)\in \Gamma_{n-1}(\xi_{n-1})$. The same being true for $y_\infty'$, we have proved $y_\infty \widehat{\mathcal R}_\infty y_\infty'$ (case ''either'' in the definition of $\widehat{\mathcal R}_\infty$).

\smallskip Case 2\,: $n=1$ and $y_\infty,y_\infty'\in \Gamma_1(x_1)$ (by Proposition \ref{proj1}, this includes the case $n=0$). 

For any $k\geq 1$ one can find $x_{0,k} \in X_0$ such that $\varphi_{k,1}(x_1)=\pi_{k,0}(y_1,y_2,\ldots,y_k,x_{0,k})$, hence $\varphi_{k,1}(x_1)=\pi_{1,k-1}\big(y_1,\pi_{k-1,0}(y_2,\ldots,y_k,x_{0,k})\big)$. By Lemma \ref{Xm}, there exists $\xi_{0,k}\in X_0$ such that $\pi_{k-1,0}(y_2,\ldots,y_k,x_{0,k})=\varphi_{k-1,0}(\xi_{0,k})$. Notice that, by Lemma \ref{X0}, this implies $\varphi_{m,0}(x_{0,k})\in C(y_2,\ldots,y_m)$ for any $m\leq k$.

Since $\{ \xi_{0,k} : k\in\bn\}\subset X_0$, which is compact Hausdorff, there is a subnet $\{ \xi_{0,k(\a)} : \a\in\ca\}$ such that $\xi_0=\lim_{\a\in\ca} \xi_{0,k(\a)}$ exists in $X_0$. On one side, $\pi_{1,0}(y_1,\xi_{0,k(\a)})=x_1$ implies $\pi_{1,0}(y_1,\xi_0)=x_1$.  On the other side, for any fixed $m\geq 1$, one has $\varphi_m(\xi_{0,k(\a)})\in C(y_2,\ldots,y_m)$ and, in the limit, $\varphi_m(\xi_0)\in C(y_2,\ldots,y_m)$.

We have proved that there exists $\xi_0\in X_0$ such that $\pi_{1,0}(y_1,\xi_0)=x_1$ and $(y_2,\ldots, y_m, \ldots)\in \Gamma_0(\xi_0)$.

\smallskip Similarly, there exists $\xi'_0\in X_0$ such that $\pi_{1,0}(y'_1,\xi'_0)=x_1$ and $(y'_2,\ldots, y'_m, \ldots)\in \Gamma_0(\xi'_0)$. We have proved $y_\infty \widehat{\mathcal R}_\infty y_\infty'$ (case ''or'' in the definition of $\widehat{\mathcal R}_\infty$).
\end{proof}

\begin{cor}\label{phi0infty}
If $\f_{\infty,0}$ is injective, all maps $\f_{\infty,n}$ are injective, hence the scheme is fully injective.
\end{cor}
\begin {proof}
Let us observe that diagram \eqref{FixedPointBis} gives the following commutative diagram for the fixed point $X_\infty$:
$$
\begin{matrix}
Y\times X_0&\stackrel{id_Y\times \f_{\infty,0}}{\hooklongrightarrow}&Y\times X_\infty
\\
\downarrow\pi_1& &\downarrow\pi_\infty\\
X_1&\stackrel{\f_{\infty,1}}{\hooklongrightarrow}&X_\infty\,.
\end{matrix}
$$
Assume that two points $x_1,x_1'$ give the same point in $X_\infty$. Then there exist $(y,x_\infty), (y',x'_\infty)\in Y\times X_\infty$ such that 
$\pi_\infty(y,x_\infty)=\pi_\infty(y',x'_\infty)$,
$x_\infty=\f_{\infty,0}(x_0)$ and $x'_\infty=\f_{\infty,0}(x'_0)$ for suitable points $x_0,x'_0\in X_0$, such that
$\pi_1(y,x_0)=x_1$, $\pi_1(y',x'_0)=x'_1$. 
Since $\f_{\infty,0}$ is injective, $x_0$ and $x'_0$ are uniquely determined, hence the commutativity of the diagram implies $x_1=x'_1$, namely $\f_{\infty,1}$ is injective.
An inductive argument concludes the proof.
\end{proof}

\section{Comparison with the standard notions of self-similar fractals}

\subsection{From fully injective similarity schemes to  Kigami-Kameyama self-similar sets}

Let us recall that according to A. Kameyama \cite{Ka},  a compact Hausdorff space $K$ is a topological self-similar set if there exist continuous maps $f_1,\dots,f_N:K\to K$ and a continuous surjection $\pi_K :\{1,\dots,N\}^\bn\to K$ such that, for all $k=1,\dots,N$ the following diagram is commutative
\begin{equation}\label{KKset}
\begin{matrix}
\{1,\dots,N\}^\bn&\stackrel{W_k}{\longrightarrow}&\{1,\dots,N\}^\bn\\
\downarrow\pi_K & & \downarrow\pi_K\\
K&\stackrel{f_k}{\longrightarrow}& K,
\end{matrix}
\end{equation}
where $W_y:\{y_\infty\}\in Y^\bn\mapsto\{y,y_\infty\}\in Y^\bn$.


\begin{rem} 
$(1)$ It follows from \cite{Wi}, Corollary 23.2, that $K$ is metrizable.

\item[$(2)$] Let us observe that an Iterated Function System, namely a complete metric space endowed with finitely many contractions, gives rise to a topological self-similar set, where the contractions give the maps $f_k$ and $K$ is the fixed point in the space of compact subspaces of the complete metric space.

\item[$(3)$] J. Kigami in \cite{Ki} introduces a particular case of topological self-similar space, where each $f_i$ is injective, and calls such $\{K;f_1,\ldots,f_N\}$ a self-similar structure.
\end{rem}


We now show that a similarity scheme gives rise to shift maps on the limit fixed point $X_\infty$ and, when $Y=\{1,\dots,N\}$, to a Kigami-Kameyama topological fractal.
Assume we have a fully injective similarity scheme $(X_0,X_1,Y)$.
 
\subsubsection{Shift endomorphisms}

Given $y\in Y$, consider the shift map $W_y : y_\infty\in Y^\bn \mapsto (y,y_\infty) \in Y^\bn$,
and the diagram
$$
\begin{matrix}
Y^\bn&\stackrel{W_y}{\longrightarrow}&Y^\bn\\
\downarrow\pi^\bn&&\downarrow\pi^\bn\\
X_\infty&& X_\infty.
\end{matrix}
$$
We now show that $W_y$ passes to the quotient.
Let $y_\infty,y'_\infty\in Y^\bn$ be such that $y_\infty\ne y'_\infty$ and $y_\infty\,\car\, y'_\infty$.
By definition, there exists $n\in\bn, x_n\in X_n$, such that $y_\infty,y'_\infty\in\Gamma_n(x_n)$. By equation \eqref{p-q-equation}, this implies that, for any $y\in Y$, $(y,y_\infty),(y,y'_\infty)\in\Gamma (\pi_{1,n}(y,x_n))$. This implies that $W_y(y_\infty)\,\car\,W_y(y'_\infty)$.

We have proved that 
\begin{prop}
For any $y\in Y$, there exists a well defined shift endomorphism $w_y:X_\infty\to X_\infty$
such that the following diagram is commutative:
\begin{equation}\label{KKshifts}
\begin{matrix}
Y^\bn&\stackrel{W_y}{\longrightarrow}&Y^\bn\\
\downarrow\pi^\bn&&\downarrow\pi^\bn\\
X_\infty&\stackrel{f_y}{\longrightarrow}& X_\infty.
\end{matrix}
\end{equation}
\end{prop}

Since the map $Y\times Y^\bn\to Y^\bn$ sending $(y,y_\infty)\mapsto \Phi_y(y_\infty)$ is surjective,  the map $\pi_\infty\stackrel{{\mathrm def}}{=}Y\times X_\infty\to X_\infty$ such that $\pi_\infty(y,x_\infty)= w_y(x_\infty)$ is surjective too.

If $y_n=(z_1,z_2,\dots,z_n)\in Y^n$, we set $W_{y_n}=W_{z_1}\circ\dots\circ W_{z_n}$ and $f_{y_n}=f_{z_1}\circ\dots\circ  f_{z_n}$, hence  diagram \eqref{KKshifts} is commutative and all maps are surjective. 
\begin{cor}\label{LimFixPtIsKKfractal}
If $Y=\{1,\dots,n\}$, the limit fixed point $X_\infty$ of the similarity scheme $(X_0,X_1,Y)$ is a Kameyama topological self-similar set.
\end{cor}

\begin{prop} 
Let $(X_0,X_1,Y,\f,\pi)$ be a fully injective similarity scheme, $y\in Y$. Then the following are equivalent
\item[$(1)$] $f_y : X_\infty \to X_\infty$ is injective, 

\item[$(2)$] for all $n\in\bn\cup\{0\}$, $\pi_{1,n}(y,\cdot) : X_n \to X_{n+1}$ is injective, 

\item[$(3)$] there is $n\in\bn\cup\{0\}$ such that $\pi_{1,n}(y,\cdot) : X_n \to X_{n+1}$ is injective.
\end{prop}
\begin{proof}
It follows from $f_y \circ \f_{\infty,n} = \f_{\infty,n+1} \circ \pi_ {1,n}(y,\cdot)$.
\end{proof}

\subsection{From Kigami self-similar structures to fully injective similarity schemes}

Given a  Kigami self-similar structure as in diagram \eqref{KKset}, we first consider the map $F$ from the subsets of $K$ given by $F(C)=\displaystyle\bigcup_{j=1,\dots,N}f_j(C)\subset K$, then call the set $S:=\{x\in K: \exists x'\in K, i,j\in \{1,\dots,N\}, i\neq j, f_i(x)=f_j(x')\}$ the essential part of the fractal $K$ and assume there exists a set $X_0\subset K$ such that
\begin{itemize}
\item[A)] $X_0$ non-empty compact set,
\item[B)] $X_0\supset S$,
\item[C)] $X_0\subset F(X_0)$,
\end{itemize}
Observe that, by injectivity of $f_j$, $X_0$ and $f_j(X_0)$ are homeomorphic.
Then set $X_1=F(X_0)$ and consider the similarity scheme 
$$
\begin{matrix}
& &\{1,\dots,N\}\times X_0
\\
& &\downarrow\pi\\
X_0&\stackrel{\f}{\hooklongrightarrow}&X_1,
\end{matrix}
$$
where $\f$ is the inclusion and $\pi$ is the map
$$
\{1,\dots,N\}\times X_0\simeq\bigsqcup_{j=1,\dots,N}f_j(X_0)
\stackrel{\pi}{\longrightarrow} \bigcup_{j=1,\dots,N}f_j(X_0)=X_1.
$$
\begin{prop}\label{Homeomorphism}
For any $k$, there exists a homeomorphism $\psi_k:F^k(X_0)\to X_k$ such that for any $k\in\bn$, the following properties hold:
\begin{itemize}
\item[$(P_k)$] $\psi_{k}(f_i(x))=\pi_{1,k-1}(i,\psi_{k-1}(x))$, $\forall x\in F^{k-1}(X_0)$,
\item[$(Q_k)$] if $x\in F^{k-1}(X_0)\subset F^{k}(X_0)$, 
then $\psi_{k}(x)=\f_{k,k-1}\circ\psi_{k-1}(x)$,
\item[$(R_k)$] $\displaystyle \forall i,j=1,\dots, n, \,x,x'\in F^k(X_0), \,
 f_i(x)=f_j(x')\Leftrightarrow \pi_{1,k}(i,\psi_k(x))=\pi_{1,k}(j,\psi_k(x')).$
\end{itemize}
\end{prop}

\begin{proof}
We first observe that, by  C), the sequence $F^k(X_0)$ is increasing. We now prove the statement by induction. This is clearly true for $k=0,1$, with $\psi_0,\psi_1$ the identity maps. Assume the statement is true for $k-1$, $k\geq2$. We first define $\psi_k$. Property $(P_k)$ forces us to set  $\psi_{k}(f_i(x))=\pi_{1,k-1}(i,\psi_{k-1}(x))$, and property $(R_{k-1})$ implies that such map is well defined and injective.
The continuity is a consequence of the following commutative diagram:
$$
\begin{matrix}
\bigsqcup_{j=1,\dots,N}f_j(F^{k-1}(X_0))&\simeq&\{1,\dots,N\}\times F^{k-1}(X_0)
&\stackrel{ id \times \psi_{k-1} }{\longrightarrow}&\{1,\dots,N\}\times X_{k-1}\\
\downarrow&&&&\downarrow\pi_{1,k-1}\\
\bigcup_{j=1,\dots,N}f_j(F^{k-1}(X_0))&=&F^k(X_0)&\stackrel{\psi_k}{\longrightarrow}& X_k.
\end{matrix}
$$
Then, being a continuous map from a compact space to a Hausdorff space, it is a homeomorphism.
\\
We now prove $(Q_k)$. For any $j=1,\dots,n$, $z\in F^{k-2}(X_0)$, we get
 \begin{align*}
  \f_{k,k-1}\circ\psi_{k-1}(f_j(z))
  &=\f_{k,k-1}\circ\pi_{1,k-2}(j,\psi_{k-2}(z))
  =\pi_{1,k-1}\circ(id \times \f_{k-1,k-2})(j,\psi_{k-2}(z))\\
 &=\pi_{1,k-1}(j,\f_{k-1,k-2}\circ\psi_{k-2}(z))
 =\pi_{1,k-1}(j,\psi_{k-1}(z))=\psi_k(f_j(z)),
 \end{align*}
 where, in the first equality, we used property $(P_{k-1})$, in the second, diagram \eqref{diagramA}, in the fourth, property $(Q_{k-1})$ and in the fifth, property $(P_{k})$.
\\
We then prove $(R_{k})$. Assume $\pi_{1,k}(i,x_k)=\pi_{1,k}(j,x'_k))$, with $(i,x_k)\ne (j,x'_k) \in \{1,\ldots,N\} \times X_k$. Then, by the construction of diagram \eqref{diagramA}, 
$\exists x_{k-1},x'_{k-1}\in X_{k-1}$ such that $\f_{k,k-1}(x_{k-1})=x_k$, $\f_{k,k-1}(x'_{k-1})=x'_k$, and 
$\pi_{1,k-1}(i,x_{k-1})=\pi_{1,k-1}(j,x'_{k-1}))$. 
Then,
\begin{align*}
\pi_{1,k}(i,x_k)
&=\pi_{1,k}\circ(id \times \f_{k,k-1})(i,x_{k-1})
=\f_{k+1,k}\circ\pi_{1,k-1}(i,x_{k-1})
=\f_{k+1,k}\circ\psi_k\circ f_i\circ\psi_{k-1}^{-1}(x_{k-1}),
\end{align*}
where, in the second equality, we used diagram \eqref{diagramA} and, in the third, property $(P_{k-1})$. Replacing $(i,x_k)$ with $(j,x'_k)$, we get 
$\f_{k+1,k}\circ\psi_k\circ f_i\circ\psi_{k-1}^{-1}(x_{k-1})=\f_{k+1,k}\circ\psi_k\circ f_j\circ\psi_{k-1}^{-1}(x'_{k-1})$. Since $\f_{k+1,k}\circ\psi_k$ is injective, we obtain
$$ f_i\circ\psi_{k-1}^{-1}(x_{k-1})=f_j\circ\psi_{k-1}^{-1}(x'_{k-1}).$$
Since
$$
f_i\circ\psi_{k-1}^{-1}(x_{k-1})=f_i\circ\psi_k^{-1}\circ\psi_k\circ\psi_{k-1}^{-1}(x_{k-1})
=f_i\circ\psi_k^{-1}\circ\f_{k,k-1}(x_{k-1})=f_i(\psi_k^{-1}(x_k)),
$$
equality $f_i(\psi_k^{-1}(x_k))=f_j(\psi_k^{-1}(x'_k))$ follows.

Now assume $f_i(z)=f_j(z')$, with $ i,j=1,\dots, n$, $z\ne z'\in F^k(X_0)$. By definition, $z,z'\in S\subset X_0$, hence, by property $(R_0)$, we get $\pi_{1,0}(i,\psi_0(z))=\pi_{1,0}(j,\psi_0(z'))$.
Now observe that, by diagram \eqref{diagramA},
$$
\pi_{1,k}(i,\psi_k(z))
=\pi_{1,k}(i,\f_{k,0}\circ\psi_0(z))
=\f_{k+1,1}\circ\pi_{1,0}(i,\psi_0(z)),
$$
and the same for $(j,z')$,
hence $\pi_{1,k}(i,\psi_k(z))=\pi_{1,k}(j,\psi_k(z'))$.
\end{proof}


Proposition \ref{Homeomorphism} immediately implies the following:

\begin{cor} \label{cor:psi}
Let $y_\infty\in\{1,\dots,N\}^\bn$. Then
\begin{itemize}
\item[$(1)$] $\psi_{n+k}(f_{(y_1,\ldots,y_n)}(x)) = \pi_{n,k}(y_1,\ldots,y_n,\psi_{k}(x))$, $x\in F^{k}(X_0)$, $n\geq1$, $k\geq0$,
\item[$(2)$]  $\psi_{p}(x) = \f_{p,k}\circ\psi_{k}(x)$, where $p\geq k\geq0$, $x\in F^{k}(X_0)\subset F^{p}(X_0)$.
\end{itemize}
\end{cor}

\begin{lem}\label{fI(x)inX0}
Let $x\in K$, $I\in\{1,\dots,N\}^k$, such that $f_{I}(x)\in X_0$. Then $x\in X_0$.
\end{lem}
\begin{proof}
(1) We first prove it for $k=1$. Indeed, $f_{i}(x)\in X_0\subset X_1=\cup_{j=1,\dots,N}f_j(X_0)$, hence there exists $1\leq j\leq N$, $x_0\in X_0$ such that $f_{i}(x)=f_{j}(x_0)$. If $i=j$, then, by injectivity, $x_0=x\in X_0$. If $i\ne j$, then by definition $x,x_0\in S\subset X_0$.
\\
(2) We then proceed by induction, assuming the property is true for $k$. If $I=(i_1,\ldots,i_{k+1})\in\{1,\dots,N\}^{k+1}$, we get $f_{(i_1,\ldots,i_k)}\circ f_{i_{k+1}}(x)\in X_0$. By the inductive hypothesis, $f_{i_{k+1}}(x)\in X_0$ and, finally,  $x\in X_0$, by (1).
\end{proof}

Let us consider two relations on $\{1,\dots,N\}^\bn$: the relation $\car$ defined in \eqref{car} and the relation $\car'$ where  $I\ \car'\, J$ iff $\pi_K(I)=\pi_K(J)$, cf. diagram \eqref{KKset}.

\begin{thm} \label{thm:KK}
Let $\{ K; f_1,\ldots,f_N\}$ be a self-similar structure, and assume there exists $X_0\subset K$  satisfying hypotheses $(A), (B),(C)$. Then, the two relations above coincide, the similarity scheme is fully injective and $K=X_\infty$.
\end{thm}
\begin{proof}
Let $I\ne J\in \{1,\dots,N\}^\bn$, with $I\,\car\, J$. Then $\exists x_n\in X_n$ such that
$$
\f_{p,n}(x_n)\in \{\pi_{p,0}(i_1,\ldots,i_p,x):x\in X_0\}\cap \{\pi_{p,0}(i_1,\ldots,i_p,x):x\in X_0\},\forall p\geq n.
$$
Then, setting $z=\psi_n^{-1}(x_n)\in F^n(X_0)$, by Corollary \ref{cor:psi} $(2)$ we get $\psi_p(z)=\f_{p,n}(x_n),\forall p\geq n$, and, by Corollary \ref{cor:psi} $(1)$,
$$
z\in f_{(i_1,\ldots,i_p)}(X_0)\cap f_{(j_1,\ldots,j_p)}(X_0),\forall p\geq n.
$$
This means that $z=\pi_K(I)=\pi_K(J)$, namely $I\,\car'\,J$.

Viceversa, assume $I\ne J\in \{1,\dots,N\}^\bn$, with $I\,\car'\,J$, namely $\pi_K(I)=\pi_K(J)=z\in K$. This implies that $z\in f_{(i_1,\ldots,i_k)}(K) \cap f_{(j_1,\ldots,j_k)}(K)$, for any $k\in\bn$. Let $n\geq1$ be the first index such that $i_{n}\ne j_{n}$, so that  $i_k=j_k$, $\forall k<n$. Since $z\in  f_{(i_1,\ldots,i_n)}(K) \cap f_{(j_1,\ldots,j_n)}(K)$,
we get $z_1,z_2\in K$ such that $z=f_{(i_1,\ldots,i_{n-1})}\circ f_{i_{n}}(z_1)=f_{(i_1,\ldots,i_{n-1})}\circ f_{j_{n}}(z_2)$, whence $f_{i_{n}}(z_1)= f_{j_{n}}(z_2)$, by the injectivity of $f_{(i_1,\ldots,i_{n-1})}$. By definition, $z_1$ and $z_2\in S\subset X_0$, from which we get $x_n\stackrel{{\rm def}}{=}\psi_n(z)=\psi_n\circ f_{(i_1,\ldots,i_n)}(z_1)=\pi_{n,0}(i_1,\ldots,i_n,z_1)=\pi_{n,0}(j_1,\ldots,j_n,z_2)\in X_{n}$.

For any $p>n$, we also have $z\in f_{(i_1,\ldots,i_p)}(K)$, namely $\exists z'\in K$ such that $z=f_{(i_1,\ldots,i_p)}(z')=f_{(i_1,\ldots,i_n)}\circ f_{(i_{n+1},\ldots,i_p)}(z')=f_{(i_1,\ldots,i_n)}(z_1)$, which implies $f_{(i_{n+1},\ldots,i_p)}(z')=z_1$. By Lemma \ref{fI(x)inX0}, $z'\in X_0$, therefore 
$$
\f_{p,n}(x_n)=\f_{p,n}\circ \psi_n(z)=\psi_p(z)=\psi_p\circ f_{(i_1,\ldots,i_p)}(z')=\pi_{p,0}(i_1,\ldots,i_p,\psi_0(z')).
$$
Therefore, $I\in \G_n(x_n)$. Since the same holds for $J$, we get $I\,\car\, J$. We proved that $\car$ and $\car'$ are the same relation, hence $K=X_\infty$.
Since $\car'$ is an equivalence relation (cf. \cite{Ka}, remark 1.3), $\car$ is an equivalence relation too.

We now prove that the scheme is fully injective. By Section \ref{section:CombDefXinfinity}, if $x_n\in X_n$, $\f_{\infty,n}(x_n)$ is the unique element given by $\pi^\bn(\Gamma_n(x_n))$.
Now, if $I\in\Gamma_n(x_n)$, $\forall p\geq n$, $\exists x_{0,p}\in X_0$ such that $\f_{p,n}(x_n) = \pi_{p,0}(i_1,\ldots,i_p,x_{0,p})$. Let us set $z_n := \psi_n^{-1}(x_n)$, so that $\psi_p(z_n) = \f_{p,n}\circ\psi_n(z_n) = \f_{p,n}(x_n) = \pi_{p,0}(i_1,\ldots,i_p,x_{0,p}) = \psi_p \circ f_{(i_1,\ldots,i_p)}(x_{0,p})$, and $z_n = f_{(i_1,\ldots,i_p)}(x_{0,p})$, which implies $z_n \in \cap_{p\geq n} f_{(i_1,\ldots,i_p)}(K)$.
Since the intersection is a singleton,  $\psi^{-1}(x_n)$, hence $x_n$, is uniquely determined by $\Gamma_n(x_n)$.
\end{proof}

\section{Examples}
The word example has two meanings here: we first show which known fractals or class of fractals can be described by a similarity scheme. Then we propose (new) similarity schemes determining if they are fully injective or not, and  describe the space $X_\infty$.

\subsection{Kigami-Kameyama topological fractals that can be described via a similarity scheme}
In this subsection, we recover examples of Kigami-Kameyama fractals for which the space $X_0$ as in Theorem \ref{thm:KK} may be constructed, so that
 there exists a similarity scheme which reproduces the given fractal.

It is worth mentioning here that  there is always a trivial solution, namely $X_0=K$. In this case the similarity scheme is already a fixed point, namely our notion is not really constructive, since the fractal is given from the outset.

Let us recall the notion of pcf self-similar structure according to Kigami, \cite{KigamiBook}, Definition 1.3.13. We define $C_K := \cap_{i=1}^N \cap_{j\neq i} (f_i(K)\cap f_j(K))$, $\cc:= \pi_K^{-1}(C_K)$, $\cp:= \cup_{n=1}^\infty \s^n(\cc)$, where $\s:(i_1,i_2,\ldots)\in \{1,\ldots,N\}^\bn\mapsto (i_2,i_3,\ldots)\in \{1,\ldots,N\}^\bn$, and finally $V_0:=\pi_K(\cp)$. Then a self-similar structure $\{K;f_1,\ldots,f_N\}$ is called post-critically finite (pcf) if $\cp$ (and hence $V_0$) is a finite set.

\begin{ex} 
If $S=\emptyset$, and $X_0$ is the set of fixed-points of $\{f_1,\ldots,f_N\}$, then the hypotheses $(A), (B), (C)$ are satisfied, the self-similar structure is pcf, we obtain a fully injective similarity scheme, and $X_\infty$ is homeomorphic to $K$ and totally disconnected.
\end{ex}
%


\begin{ex} 
If $S\neq\emptyset$, and the self-similar structure is pcf, we choose $X_0=V_0$. Then the hypotheses $(A), (B), (C)$ are satisfied, we obtain a fully injective similarity scheme, and $X_\infty$ is homeomorphic to $K$. 
If $V_0$ is not finite, we may put $X_0=\ov{V_0}$, however such a similarity scheme may be trivial, namely equality $X_0=X_\infty$ may hold. 

\end{ex}
%


\begin{ex} 
If $S\neq\emptyset$, and $S\subset \Phi(S)$, we choose $X_0=S$. Then the hypotheses $(A), (B), (C)$ are satisfied, we obtain a fully injective similarity scheme, and $X_\infty$ is homeomorphic to $K$.

This family of examples contains non pcf self-similar structures like the Sierpinski carpet. Let us recall that the Sierpinski carpet is the fixed point of an iterated function system, where $\Omega$ is given by $\br^2$ with the standard metric and one has 8 contraction maps with scaling parameter 1/3, 4 of which have fixed points in the vertices of a unit square, and 4 of which have fixed points in the middle points of the four sides of the square.
According to the definition of $S$ in the section above, $S$ for the Sierpinski carpet is given by the boundary of the unit square, therefore $\Phi(S)$, being given by the union of 8 rescaled copies of $S$, actually contains $S$.
\end{ex}

\subsection{Further examples: new similarity schemes}

\begin{ex} [The diagonal similarity scheme] 
Let $Y$ be any Hausdorff compact set, and set $X_0=Y$, $X_1=Y\times Y/\sim$ where the equivalence relation is the folding along the diagonal, namely we assume $(a,b)\sim(b,a)$ for any $a,b\in Y$. We also define $\f:Y\to X_1$ by $y\mapsto (y,y)$.
We observe that this similarity scheme is fully injective, namely the map $\f_{\infty,0}$ from $X_0$ to the limit fixed point $X_\infty$ is injective. Indeed, $X_\infty$ is a quotient of $Y^\bn$ with respect to an infinite set of equivalence relations: $(y_1,y_2,y_\infty)\sim(y_2,y_1,y_\infty)$,
$(y_1,y_2,y_3,y_\infty)\sim (y_1,y_3,y_2,y_\infty)$, $(y_1,y_2,y_3,y_4,y_\infty)\sim
(y_1,y_2,y_4,y_3,y_\infty)$ and so on. In particular, $(y,y,y,\dots)$ has no other elements in its equivalence class, and $\f_{0,\infty}(y)=(y,y,y,\dots)$, namely it is injective.
By Corollary \ref{phi0infty} the scheme is fully injective.

Let us observe that when $Y=X_0$ is a finite set the diagonal similarity scheme reproduces topologically some known IFS's: if $|Y|=2$, our limit fractal is homeomorphic to the  interval $[0,1]$ obtained as a fixed point in $\br$ of the contractions centered in the points $0$ and $1$ with scaling parameter 1/2. If $|Y|>2$, we obtain a gasket inscribed in the equilateral $n$-simplex, described as the fixed point in $\br^{n-1}$ of $n$-contractions centered in the vertices of the simplex with scaling parameter 1/2. In particular, for $n=3$, we get the Sierpinki gasket.
\end{ex}


%

\begin{rem}
The theory described in this paper concerns fully injective similarity schemes. Indeed, strictly speaking, the notion of fixed point $(Z,\f_Z)$ of a similarity scheme requires  $\f_Z$ to be injective. However, it is not difficul to see that, extending the category to pairs where injectivity is not required, the notion of fixed point still has a meaning. In some cases limit points of the similarity scheme are still fixed points, as the following examples show. Moreover, uniqueness of the fixed point may fail, as the last example shows.
\end{rem}

\begin{ex} [A similarity scheme which is not fully injective] 
Let $Y=\{0,1\}$, $X_0=\{a,b\}$, $X_1=Y\times X_0/\sim$ where $(0,a)\sim(1,a)$, $\f(a)=(0,a)$, $\f(b)=(0,b)$. 
Clearly, $Y^n\times X_0$ consists of the elements $\{(r_n,a), (s_n,b), r_n,s_n\in \{0,1\}^n\}$ and the equivalence classes associated with the projection $\pi_{n,1}$ consists of the set  $\{(r_n,a),  r_n,\in \{0,1\}^n\}$ and of the singletons $(s_n,b)$, $s_n\in \{0,1\}^n\}$. Denoting by $a_n$ the equivalence class $\{(r_n,a),  r_n,\in \{0,1\}^n\}$, it turns out that $X_n=\{a_n,(s_n,b):s_n\in \{0,1\}^n\}$ and $\f_{n,n+1}(a_n)=a_{n+1}$, $\f_{n,n+1}(s_n,b)=((s_n,0),b)$. Then, if $y_\infty\in \{0,1\}^\bn$, $c(\rho_n(y_\infty))=\{a_n,(\rho_n(y_\infty),b)\}$, namely $\f_{n,0}(a)=a_n\in c(\rho_n(y_\infty))$ for any $n\in \bn, y_\infty\in Y^\bn$, or, equivalently, $\Gamma(a)=Y^\bn$. As a consequence, $X_\infty$ is a singleton. We observe that, even though the map $\f_{\infty,0}:X_0\to X_\infty$ is not injective, $(X_\infty,\f_{\infty,0})$ is still a fixed point of the endo-functor
 $\Phi$.
\end{ex}


The following example is an adaptation of Example 1 of Section 3 in \cite{Sece12}.
 
\begin{ex} [A similarity scheme with more than one fixed-point] \label{nonunique}  
Let $Y=\{0,1\}$, $X_0=\{a,b,c\}$, $X_1=Y\times X_0/\sim$ where $(0,a)\sim(1,a)$, $(0,c)\sim(1,c)$, $\f(x_0)=[(0,x_0)]_\sim$. 
Reasoning as in the previous example, $X_\infty$ is a singleton. Moreover, if $Z:=\{a,c\}$, and $\f_Z:X_0\to Z$ is defined as $\f_Z(a)=\f_Z(b)=a$, $\f_Z(c)=c$, one can prove that $(Z,\f_Z)$ is another fixed point of the functor $\Phi$.
\end{ex}



\begin{ack}This work was partially supported by the ERC Advanced Grant 669240 QUEST \textquotedblleft Quantum Algebraic Structures and Models\textquotedblright\,.
We also acknowledge the support of GNAMPA-INdAM for financing the stay of J-LS in Rome in October 2023.
D. G. and T. I. acknowledge the MUR Excellence Department Project awarded to the Department of Mathematics, University of Rome Tor Vergata, MatMod@TOV (CUP E83C23000330006)
 and the University of Rome Tor Vergata funding scheme 
``OAQM''  CUP E83C22001800005.
\end{ack}

The results contained in this paper were presented in various international conferences, by D.G. in the conference ``Analysis on fractals and networks, and applications'', 18 - 22 March, 2024, Marseille, France; By T.I. in the conference ``Fractals, quantum graphs, and their applications in pure and applied sciences'', 25 - 27 March 2024, Politecnico di Milano, Italy; and by J.-L. S. in the conference ``Noncommutative Geometry and Applications'', 24 - 28 June 2024, Cortona, Italy

\end{document}